\documentclass[reqno,10pt]{amsart}

\usepackage[all]{xy}

\theoremstyle{plain}
  \newtheorem{theorem}{Theorem}
  \newtheorem*{theo_genp}{Theorem~\ref{T:gen}'}
  \newtheorem*{theo_torp}{Theorem~\ref{T:tor}'}
  
  \newtheorem{lemma}{Lemma}
  \newtheorem{proposition}{Proposition}
  
  \newtheorem*{conjecture*}{Conjecture}
  \newtheorem{conjecture}{Conjecture}

\theoremstyle{definition}
  \newtheorem{definition}{Definition}
  
  \newtheorem*{convention*}{Convention}
  
\theoremstyle{remark}
  \newtheorem{remark}{Remark}
  
  \newtheorem{example}{Example}

\newcommand{\LL}{\textsf{L}}
\newcommand{\MS}{\textsf{MS}}
\newcommand{\NCT}{\textsf{NCT}}
\newcommand{\EG}{\textsf{EG}}
\newcommand{\SEG}{\textsf{SEG}}

\newcommand{\sm}{\text{\rm sm}}

\newcommand{\la}{\text{\rm la}}
\newcommand{\an}{\text{\rm an}}

\newcommand\R{\mathbb R}
\newcommand\QX{\mathfrak I}
\newcommand\PX{\mathfrak P}
\newcommand\RX{\mathfrak R}
\newcommand\SX{\Sigma}
\DeclareMathOperator{\ind}{ind}

\DeclareMathOperator{\eend}{end}
\DeclareMathOperator{\supp}{supp}
\DeclareMathOperator{\Fix}{Fix}
\DeclareMathOperator{\img}{img}

\DeclareMathOperator{\ev}{ev}

\DeclareMathOperator{\length}{length}
\DeclareMathOperator{\Tor}{Tor}
\DeclareMathOperator{\dist}{dist}

\DeclareMathOperator{\sign}{sign}

\DeclareMathOperator{\IND}{IND}
\DeclareMathOperator{\Vol}{Vol}
\DeclareMathOperator{\vvol}{vol}
\DeclareMathOperator{\Maps}{Maps}

\DeclareMathOperator{\Int}{Int}
\DeclareMathOperator{\grad}{grad}

\DeclareMathOperator{\tr}{tr}
\DeclareMathOperator{\diver}{div}
\DeclareMathOperator{\Gr}{Gr}
\DeclareMathOperator{\id}{id}
\DeclareMathOperator{\End}{End}
\DeclareMathOperator{\cs}{cs}
\DeclareMathOperator{\ec}{e}

\newcommand\itemref[1]{(\ref{#1})}

\begin{document}

\title[Dynamics, Laplace transform and spectral geometry]
      {Dynamics, Laplace transform and spectral geometry$^*$}

\author{Dan Burghelea}

\address{Dan Burghelea,
         Dept. of Mathematics, 
         The Ohio State University, 
         231 West Avenue, Columbus, OH 43210, USA.}

\email{burghele@math.ohio-state.edu}

\author{Stefan Haller}

\address{Stefan Haller,
         Department of Mathematics, University of Vienna,
         Nordbergstra{\ss}e 15, A-1090, Vienna, Austria.}

\email{stefan.haller@univie.ac.at}

\thanks{$^*$This paper is a new version of work which has circulated as preprint
        \cite{BH04} under the name ``Laplace transform, dynamics and spectral geometry.''
        The mathematical perspective of this version is however very different.}
\thanks{Part of this work was done while the second author enjoyed the
        warm hospitality of The Ohio State University. The second author was
        partially supported by the \emph{Fonds zur F\"orderung der 
        wis\-sen\-schaft\-lichen Forschung} (Austrian Science Fund),
        project number {\tt P14195-MAT}. Part of this work was carried out
        when the first author author enjoyed the hospitality of IHES in 
        Bures sur Yvette and second author of the Max Planck Institute for Mathematics in Bonn}



\date{\today}

\begin{abstract}
We consider a vector field $X$ on a closed manifold which admits a 
Lyapunov one form. We assume $X$ has Morse type
zeros, satisfies the Morse--Smale transversality condition and
has non-degenerate closed trajectories only.
For a closed one form $\eta$, considered as flat connection on the trivial
line bundle, the differential of the Morse complex formally associated 
to $X$ and $\eta$ is given by infinite series.
We introduce the exponential growth condition and show that it guarantees
that these series converge absolutely for a non-trivial set of $\eta$.
Moreover the exponential growth condition guarantees that we have an
integration homomorphism from the deRham complex 
to the Morse complex. We show that the integration induces an isomorphism
in cohomology for generic $\eta$. Moreover, we define a complex valued Ray--Singer
kind of torsion of the integration homomorphism, and compute it
in terms of zeta functions of closed trajectories of $X$. 
Finally, we show that the set of
vector fields satisfying the exponential growth condition is
$C^0$--dense.
\end{abstract}

\maketitle
\setcounter{tocdepth}{1}
\tableofcontents

\section{Introduction}\label{S:intro}

Let $M$ be a closed smooth manifold. We consider a vector field $X$
which admits a Lyapunov form, see Definition~\ref{D:Lyap}. We assume $X$
has Morse type zeros, satisfies Morse--Smale transversality and has
non-degenerate closed trajectories only. 
These assumptions imply that the number of instantons as well as
the number of closed trajectories in a fixed homotopy class are finite.
Moreover, we assume that $X$ satisfies the exponential growth condition, 
a condition on the growth of the volume of the unstable manifolds of $X$,
see Definition~\ref{D:0} below.
Using a theorem of Pajitnov we show that the set of vector fields with these
properties is $C^0$--dense, see Theorem~\ref{T:gen}.

Let $\eta\in\Omega^1(M;\mathbb C)$ be a closed one form, and consider it
as a flat connection on the trivial bundle $M\times\mathbb C\to M$. 
Using the zeros and instantons of $X$ one might try to associate
a More type complex to $X$ and $\eta$.
Since the number of instantons between zeros of $X$ is in general infinite,
the differential in such a complex
is given by infinite series. The exponential growth condition guarantees
that this series converges absolutely for a non-trivial set of closed one
forms $\eta$. For these $\eta$ we thus have a Morse complex
$C^*_\eta(X;\mathbb C)$, see section~\ref{SS:Mc}, which, as a `function' of $\eta$, 
can be considered as the `Laplace transform'
of the Novikov complex. The exponential growth condition also guarantees
that we have an integration homomorphism 
$\Int_\eta:\Omega^*_\eta(M;\mathbb C)\to C^*_\eta(X;\mathbb X)$,
where $\Omega^*_\eta(M;\mathbb C)$ denotes the deRham complex
associated with the flat connection $\eta$.
It turns out that this integration homomorphism induces an isomorphism
in cohomology, for generic $\eta$. These results are the contents of 
Theorem~\ref{T:int} and Proposition~\ref{P:expR}.

For those $\eta$ for which $\Int_\eta$ induces an isomorphism in cohomology
we define the (relative) torsion of $\Int_\eta$ with the help of zeta 
regularized determinants of Laplacians in the spirit of Ray--Singer. 
Our torsion however is based on non-positive Laplacians, is complex
valued, and depends holomorphically on $\eta$. While the definition requires the choice of a Riemannian
metric on $M$ we add an appropriate correction term which 
causes our torsion to be independent of this choice, see Proposition~\ref{P:Tint}.
Combining results of Hutchings--Lee, Pajitnov and Bismut--Zhang
we show that the torsion of $\Int_\eta$ coincides with the 
`Laplace transform' of the counting function for closed trajectories of $X$,
see Theorem~\ref{T:tor}. Implicitly, the set of closed one forms $\eta$
for which the Laplace transform of the counting function for closed
trajectories converges absolutely is non-trivial, providing an (exponential) 
estimate on the growth of the number of closed trajectories in each homology
class, as the class varies in $H_1(M;\mathbb Z)/\Tor(H_1(M;\mathbb Z))$.
Moreover, the torsion of $\Int_\eta$ provides an analytic continuation
of this Laplace transform, considered as a function on the space of closed
one forms, beyond the set of $\eta$ for which it is naturally defined.

The rest of the paper is organized as follows. The remaining part 
of section~\ref{S:intro} contains a thorough explanation of the main
results including all necessary definitions. The proofs are
postponed to sections~\ref{S:exp} through \ref{S:Ttor} and two appendices.

\subsection{Morse--Smale vector fields}\label{SS:MS}

Let $X$ be a smooth vector field on a smooth manifold $M$ of dimension $n$. A point
$x\in M$ is called a \emph{rest point} or a \emph{zero} if $X(x)=0$. 
The collection of these points  will be denoted by $\mathcal X:=\{x\in M\mid X(x)=0\}$.

Recall that a rest point $x\in\mathcal X$ is said to be of \emph{Morse type}  
if there exist coordinates $(x_1,\dotsc,x_n)$ centered at $x$ so that 
\begin{equation}\label{E:4}
X=\sum_{i\leq q}x_i\frac\partial{\partial x_i}-\sum_{i>q}x_i\frac\partial{\partial x_i}.
\end{equation}
The integer $q$ is called the \emph{Morse index} of $x$ and denoted by  $\ind(x)$. A 
rest point of Morse type is non-degenerate and its Hopf index is $(-1)^{n-q}$.
The Morse index  is independent of the chosen coordinates $(x_1,\dotsc,x_n)$.
Denote by $\mathcal X_q$ the set of rest points of Morse index $q$.
Clearly, $\mathcal X=\bigsqcup_q\mathcal X_q$.

\begin{convention*}
Unless explicitly mentioned all vector fields in this paper are assumed 
to have all rest points of Morse type, hence isolated.
\end{convention*}

For $x\in\mathcal X$, the \emph{stable} resp.\ 
\emph{unstable set} is defined by
$$
D^\pm_x:=\{y\mid\lim_{t\to\pm\infty}\Psi_t(y)=x\}
$$
where $\Psi_t:M\to M$ denotes the flow of $X$ at time $t$. The stable and unstable sets are
images of injective smooth immersions $i^\pm_x:W^\pm_x\to M$. 
The manifold $W^-_x$ resp.\ $W^+_x$ is diffeomorphic to $\mathbb R^{\ind(x)}$
resp.\ $\mathbb R^{n-\ind(x)}$.

\begin{definition}[Morse--Smale property, \MS]\label{D:MS}
The vector field $X$ is said to satisfy the \emph{Morse--Smale property}, \MS\ for short,
if the maps $i^-_x$ and $i^+_y$ are transversal, for all $x,y\in\mathcal X$.
\end{definition}

If the vector field $X$ satisfies \MS, 
and $x\neq y\in\mathcal X$, then the set $D^-_x\cap D^+_y$,  
is the image of an injective immersion of a smooth manifold 
$\mathcal M(x,y)$ of dimension $\ind(x)-\ind(y)$. Moreover, $\mathcal
M(x,y)$ is equipped with a free and proper $\mathbb R$--action. 
The quotient is a smooth manifold $\mathcal T(x,y)$ of dimension $\ind(x)-\ind(y)-1$,
called the manifold of \emph{trajectories} from $x$ to $y$. 
Recall that a collection $\mathcal O=\{\mathcal O_x\}_{x\in\mathcal X}$ of
orientations of the unstable manifolds, $\mathcal O_x$ being an
orientation of $W_x^-$, provides (coherent) orientations on
$\mathcal M(x,y)$ and $\mathcal T(x,y)$.
If $\ind(x)-\ind(y)=1$ then $\mathcal T(x,y)$ is zero dimensional and its elements 
are isolated trajectories called \emph{instantons.} The orientations
$\mathcal O$ provide a sign $\epsilon^{\mathcal O}(\sigma)\in\{\pm1\}$ 
for every instantons $\sigma\in\mathcal T(x,y)$.

\subsection{Closed trajectories}\label{SS:clotra}

Recall that a \emph{parameterized closed trajectory} is a 
pair $({\theta},T)$ consisting of a non-constant smooth curve
$\theta:\mathbb R\to M$ and a real number $T$ such that 
$\theta'(t)=X(\theta(t))$ and $\theta(t+T)=\theta(t)$ hold for all
$t\in\mathbb R$.  
A \emph{closed trajectory} is an equivalence class $\sigma$ of 
parameterized closed trajectories, where two parametrized closed
trajectories $(\theta_1,T_1)$ and $(\theta_2,T_2)$ are equivalent if 
there exists $a\in\mathbb R$ such that $T_1=T_2$ and 
$\theta_1(t)=\theta_2(t+a)$, for all $t\in\mathbb R$.
Recall that the \emph{period} $p(\sigma)$ of a closed trajectory $\sigma$
is the largest integer $p$ such that for some (and hence
every) representative $(\theta,T)$ of $\sigma$ the map $\theta:\mathbb
R/T\mathbb Z=S^1\to M$ factors through a map $S^1\to S^1$ of degree $p$.
Also note that every closed trajectory gives rise to a homotopy class
in $[S^1,M]$.

Suppose $(\theta,T)$ is a parametrized closed trajectory and $t_0\in\mathbb R$.
Then the differential of the flow 
$T_{\theta(t_0)}\Psi_T:T_{\theta(t_0)}M\to T_{\theta(t_0)}M$ 
fixes $X(\theta(t_0))$ and hence descends to a linear isomorphism
$A_{\theta(t_0)}$ on the normal space to the trajectory 
$T_{\theta(t_0)}M/\langle X(\theta(t_0))\rangle$, called the \emph{return map.}
Note that the conjugacy class of
$A_{\theta(t_0)}$ only depends on the closed trajectory represented by
$(\theta,T)$. Recall that a
closed trajectory is called \emph{non-degenerate} if $1$ is not an eigen
value of the return map. Every non-degenerate closed trajectory
$\sigma$ has a sign $\epsilon(\sigma)\in\{\pm1\}$
defined by $\epsilon(\sigma):=\sign\det(\id-A_{\theta(t_0)})$
where $t_0\in\mathbb R$ and $(\theta,T)$ is any representative of $\sigma$.

\begin{definition}[Non-degenerate closed trajectories, \NCT]\label{D:NCT} 
A vector field is said to satisfies the \emph{non-degenerate closed trajectories
property}, \NCT\ for short, if all of its closed trajectories are non-degenerate.
\end{definition}

\subsection{Lyapunov forms}\label{SS:Lyap}

The existence of a Lyapunov form for a vector field has several important
implications: it implies finiteness properties for
the number of instantons and closed trajectories, see
Propositions~\ref{P:Nov} and \ref{P:Hut} below; and it permits
to complete the unstable manifolds to manifolds with corners, see
Theorem~\ref{T:6} in section~\ref{SS:TUS}.

\begin{definition}[Lyapunov property, \LL]\label{D:Lyap}
A closed one form $\omega\in\Omega^1(M;\mathbb R)$ for which
$\omega(X)<0$ on $M\setminus\mathcal X$ is called \emph{Lyapunov form} for $X$.
We say a vector field satisfies the \emph{Lyapunov property}, 
\LL\ for short, if it admits Lyapunov forms.
A cohomology class in $H^1(M;\mathbb R)$ is called \emph{Lyapunov cohomology
class for $X$} if it can be represented by a Lyapunov form for $X$.
\end{definition}

The Kupka--Smale theorem \cite{K63, S63, P67} immediately implies

\begin{proposition}\label{P:KS}
Suppose $X$ satisfies \LL, and let $r\geq1$.
Then, in every $C^r$--neigh\-bor\-hood of $X$, there exists a vector field
which coincides with $X$ in a neighborhood of $\mathcal X$, and which
satisfies \LL, \MS\ and \NCT.
\end{proposition}

In appendix~\ref{S:cano} we will prove

\begin{proposition}\label{P:cano}
Every Lyapunov cohomology class for $X$ can be represented by a 
closed one form $\omega$, so that there exists a Riemannian metric $g$
with $\omega=-g(X,\cdot)$. Moreover,
one can choose $\omega$ and $g$ to have standard form in a neighborhood of
$\mathcal X$, i.e.\ locally around every zero of $X$, with respect to the 
coordinates $(x_1,\dotsc,x_n)$ in which $X$ has the form \eqref{E:4},
we have $\omega=-\sum_{i\leq q}x_idx^i+\sum_{i>q}x_idx^i$
and $g=\sum_i(dx^i)^2$. 
\end{proposition}

For the structure of the set of Lyapunov cohomology classes we obviously have

\begin{proposition}\label{P:cone}
The set of Lyapunov cohomology classes for $X$ constitutes an open convex
cone in $H^1(M;\mathbb R)$.
Consequently we have:
If $X$ satisfies \LL, then it admits a Lyapunov class contained in the
image of $H^1(M;\mathbb Z)\to H^1(M;\mathbb R)$.
If $X$ satisfies \LL, then it admits a Lyapunov class $\xi$ such that 
$\xi:H_1(M;\mathbb Z)/\Tor(H_1(M;\mathbb Z))\to\mathbb R$ is injective.
If $0\in H^1(M;\mathbb R)$ is a Lyapunov class for $X$ then
every cohomology class in $H^1(M;\mathbb R)$ is Lyapunov for $X$.
\end{proposition}

The importance of Lyapunov forms stems from the following two results.
Both propositions are a consequence of the fact that the energy 
of an integral curve $\gamma$ of $X$ satisfies 
$E_g(\gamma)=-\omega(\gamma)$ where $g$ and $\omega$ are as in
Proposition~\ref{P:cano}.

\begin{proposition}[Novikov \cite{N93}]\label{P:Nov}
Suppose $X$ satisfies \MS, let $\omega$ be a Lyapunov form for $X$,
let $x,y\in\mathcal X$ with $\ind(x)-\ind(y)=1$, and let $K\in\mathbb R$.
Then the number of instantons $\sigma$ from $x$ to $y$ which satisfy
$-\omega(\sigma)\leq K$ is finite.
\end{proposition}

\begin{proposition}[Fried \cite{F87}, Hutchings--Lee \cite{HL99}]\label{P:Hut}
Suppose $X$ satisfies \MS\ and \NCT, let $\omega$ be Lyapunov for $X$,
and let $K\in\mathbb R$.
Then the number of closed trajectories $\sigma$ which satisfy
$-\omega(\sigma)\leq K$ is finite.
\end{proposition}

\subsection{Counting functions and their Laplace transform}\label{SS:cf}

Let us introduce the notation $\mathcal Z^1(M;\mathbb
C):=\{\eta\in\Omega^1(M;\mathbb C)\mid d\eta=0\}$. Similarly, we will write
$\mathcal Z^1(M;\mathbb R)$ for the set of real valued closed one forms.
For a homotopy class $\gamma$ of paths joining two (rest) points in $M$ 
and $\eta\in\mathcal Z^1(M;\mathbb C)$ we
will write $\eta(\gamma):=\int_\gamma\eta$.

For a vector field $X$ which satisfies \LL\ and \MS, 
and two zeros $x,y\in\mathcal X$ with $\ind(x)-\ind(y)=1$, we define the
\emph{counting function of instantons from $x$ to $y$} by
$$
\mathbb I_{x,y}=\mathbb I_{x,y}^{X,\mathcal O}:\mathcal P_{x,y}\to\mathbb Z,
\qquad
\mathbb I_{x,y}(\gamma):=\sum_{\sigma\in\gamma}\epsilon^{\mathcal O}(\sigma).
$$
Here $\mathcal P_{x,y}$ denotes the space of homotopy classes of paths from
$x$ to $y$, and the sum is over all instantons $\sigma$ in the homotopy class
$\gamma\in\mathcal P_{x,y}$. Note that these sums are finite in view of
Proposition~\ref{P:Nov}. For notational simplicity we set $\mathbb I_{x,y}:=0$
whenever $\ind(x)-\ind(y)\neq1$.

Consider the `Laplace transform' of $\mathbb I_{x,y}$,
\begin{equation}\label{E:LIdef}
L(\mathbb I_{x,y}):\QX_{x,y}\to\mathbb C,\qquad
L(\mathbb I_{x,y})(\eta)
:=\sum_{\gamma\in\mathcal P_{x,y}}\mathbb I_{x,y}(\gamma)e^{\eta(\gamma)}
\end{equation}
where $\QX_{x,y}=\QX^X_{x,y}\subseteq\mathcal Z^1(M;\mathbb C)$ 
denotes the subset of closed one forms $\eta$ for which 
this sum converges absolutely. Moreover, set 
$\QX:=\bigcap_{x,y\in\mathcal X}\QX_{x,y}$, and let $\mathring\QX_{x,y}$ resp.\ 
$\mathring\QX=\bigcap_{x,y\in\mathcal X}\mathring\QX_{x,y}$ denote the interior of
$\QX_{x,y}$ resp.\ $\QX$ in $\mathcal Z^1(M;\mathbb C)$
equipped with the $C^\infty$--topology.

Classically \cite{W46} the Laplace transform is a partially defined 
holomorphic function $z\mapsto\int_{\mathbb R}e^{-z\lambda}d\mu(\lambda)$, associated to a complex
valued measure $\mu$ on the real line with support bounded from below.
The Laplace transform has an abscissa of absolute convergence $\rho\leq\infty$
and will converge absolutely for $\Re(z)>\rho$.
If the measure has discrete support this specializes to Dirichlet series,
$z\mapsto\sum_ia_ie^{-z\lambda_i}$.

One easily derives the following proposition which
summarizes some basic properties of $L(\mathbb I_{x,y}):\QX_{x,y}\to\mathbb C$ 
analogous to basic properties of classical Laplace transforms \cite{W46}.
The convexity follows from H\"older's inequality.

\begin{proposition}\label{P:basic:I}
The set $\QX_{x,y}$ (and hence $\mathring\QX_{x,y}$)
is convex and we have $\QX_{x,y}+\omega\subseteq\QX_{x,y}$ for all
$\omega\in\mathcal Z^1(M;\mathbb C)$ with $\Re(\omega)\leq0$.
Moreover, $\QX_{x,y}$ and \eqref{E:LIdef} are gauge invariant, i.e.\
for $h\in C^\infty(M;\mathbb C)$ and $\eta\in\QX_{x,y}$ we have
$\QX_{x,y}+dh\subseteq\QX_{x,y}$ and
$$
L(\mathbb I_{x,y})(\eta+dh)=L(\mathbb I_{x,y})(\eta)e^{h(y)-h(x)}.
$$
The restriction $L(\mathbb I_{x,y}):\mathring\QX_{x,y}\to\mathbb C$ is
holomorphic.\footnote{For a definition of holomorphicity in infinite
dimensions see \cite{HP74}.}
If $\omega$ is Lyapunov for $X$ then $\QX_{x,y}+\omega\subseteq\mathring\QX_{x,y}$, 
and for all $\eta\in\QX_{x,y}$
\begin{equation}\label{E:cont:I}
\lim_{t\to0^+}L(\mathbb I_{x,y})(\eta+t\omega)=L(\mathbb I_{x,y})(\eta).
\end{equation}
Particularly, $\mathring\QX_{x,y}\subseteq\QX_{x,y}$ is dense,
and the function $L(\mathbb I_{x,y}):\QX_{x,y}\to\mathbb C$ is completely
determined by its restriction to $\mathring\QX_{x,y}$.
\end{proposition}

\begin{remark}
In view of the gauge invariance $L(\mathbb I_{x,y})$ can be regarded as a
partially defined holomorphic function on the finite dimensional vector
space $H^1(M;\mathbb C)\times\mathbb C$.
\end{remark}

For a vector field $X$ which satisfies \LL, \MS\ and \NCT\
we define its \emph{counting function of closed trajectories} by
$$
\mathbb P=\mathbb P^X:[S^1,M]\to\mathbb Q,
\qquad
\mathbb P(\gamma):=\sum_{\sigma\in\gamma}
\frac{\epsilon(\sigma)}{p(\sigma)}.
$$
Here $[S^1,M]$ denotes the space of homotopy classes of maps $S^1\to M$,
and the sum is over all closed trajectories $\sigma$ in the homotopy class
$\gamma\in[S^1,M]$. Note that these sums are finite in view of
Proposition~\ref{P:Hut}. Moreover, define
$$
h_*\mathbb P:H_1(M;\mathbb Z)/\Tor(H_1(M;\mathbb Z))\to\mathbb Q,
\qquad
(h_*\mathbb P)(a):=\sum_{h(\gamma)=a}\mathbb P(\gamma)
$$
where $h:[S^1,M]\to H_1(M;\mathbb Z)/\Tor(H_1(M;\mathbb Z))$, and the sum is over all
$\gamma\in[S^1,M]$ for which $h(\gamma)=a$. Note that these are finite sums
in view of Proposition~\ref{P:Hut}.

Consider the `Laplace transform' of $h_*\mathbb P$,
\begin{equation}\label{E:LPdef}
L(h_*\mathbb P):\PX\to\mathbb C,\qquad
L(h_*\mathbb P)(\eta)
:=\sum_{a\in H_1(M;\mathbb Z)/\Tor(H_1(M;\mathbb Z))}(h_*\mathbb P)(a)e^{\eta(a)}
\end{equation}
where $\PX=\PX^X\subseteq\mathcal Z^1(M;\mathbb C)$ 
denotes the subset of closed one forms $\eta$ for which 
this sum converges absolutely.\footnote{We will see that $L(h_*\mathbb P)(\eta)$ 
converges absolutely in some interesting cases, see Theorem~\ref{T:tor}
below. However, our arguments do not suffice to prove (absolute) convergence
of $L(\mathbb P)(\eta):=\sum_{\gamma\in[S^1,M]}\mathbb
P(\gamma)e^{\eta(\gamma)}$. 
Of course $L(h_*\mathbb P)(\eta)=L(\mathbb P)(\eta)$, provided the latter
converges absolutely.} Let $\mathring\PX$ denote the interior of
$\PX$ in $\mathcal Z^1(M;\mathbb C)$ equipped with the $C^\infty$--topology.
Analogously to Proposition~\ref{P:basic:I} we have

\begin{proposition}\label{P:basic:P}
The set $\PX$ (and hence $\mathring\PX$) is convex and we have
$\PX+\omega\subseteq\PX$ for all $\omega\in\mathcal Z^1(M;\mathbb C)$ with 
$\Re(\omega)\leq0$. Moreover, $\PX$ and \eqref{E:LPdef} are gauge invariant,
i.e.\ for $h\in C^\infty(M;\mathbb C)$ and $\eta\in\PX$ we have
$$
L(h_*\mathbb P)(\eta+dh)=L(h_*\mathbb P)(\eta).
$$
The restriction $L(h_*\mathbb P):\mathring\PX\to\mathbb C$ is holomorphic.
If $\omega$ is Lyapunov for $X$ then $\PX+\omega\subseteq\mathring\PX$, 
and for all $\eta\in\PX$
\begin{equation}\label{E:cont:P}
\lim_{t\to0^+}L(h_*\mathbb P)(\eta+t\omega)=L(h_*\mathbb P)(\eta).
\end{equation}
Particularly, $\mathring\PX\subseteq\PX$ is dense,
and the function $L(h_*\mathbb P):\PX\to\mathbb C$ is completely determined
by its restriction to $\mathring\PX$.
\end{proposition}

\begin{remark}
In view of the gauge invariance $L(h_*\mathbb P)$ can be regarded as
a partially defined holomorphic function on the finite dimensional vector
space $H^1(M;\mathbb C)$.
\end{remark}

For $x\in\mathcal X$ let $L^1(W_x^-)$
denote the space of absolutely integrable functions $W_x^-\to\mathbb C$ with 
respect to the measure induced from the Riemannian metric $(i_x^-)^*g$,
where $g$ is a Riemannian metric on $M$.
The space $L^1(W_x^-)$ does not depend on $g$.
For a closed one form $\eta\in\mathcal Z^1(M;\mathbb C)$ 
let $h^\eta_x:W_x^-\to\mathbb C$ denote the unique smooth function which
satisfies $h^\eta_x(x)=0$ and $dh_x^\eta=(i_x^-)^*\eta$.
For $x\in\mathcal X$ define 
$$
\RX_x=\RX^X_x:=
\bigl\{\eta\in\mathcal Z^1(M;\mathbb C)\bigm| e^{h^\eta_x}\in
L^1(W_x^-)\bigr\},
$$
and set $\RX:=\bigcap_{x\in\mathcal X}\RX_x$.
Moreover, let $\mathring\RX_x$ resp.\ $\mathring\RX=\bigcap_{x\in\mathcal X}\mathring\RX_x$ 
denote the interior of $\RX_x$ 
resp.\ $\RX$ in $\mathcal Z^1(M;\mathbb C)$ equipped with the $C^\infty$--topology.

For $\alpha\in\Omega^*(M;\mathbb C)$ 
consider the `Laplace transform' of
$(i_x^-)^*\alpha\in\Omega^*(W_x^-;\mathbb C)$,
\begin{equation}\label{E:LRdef}
L((i_x^-)^*\alpha):\RX_x\to\mathbb C,\qquad
L((i_x^-)^*\alpha)(\eta):=\int_{W_x^-}e^{h_x^\eta}\cdot(i_x^-)^*\alpha.
\end{equation}
Note that these integrals converge absolutely for $\eta\in \RX_x$.
Analogously to Propositions~\ref{P:basic:I} and \ref{P:basic:P} we have

\begin{proposition}\label{P:basic:R}
The set $\RX_x$ (and hence $\mathring\RX_x$) is convex and we have $\RX_x+\omega\subseteq\RX_x$ for all
$\omega\in\mathcal Z^1(M;\mathbb C)$ with $\Re(\omega)\leq0$.
Moreover, $\RX_x$ and \eqref{E:LRdef} are gauge invariant, i.e.\
for $h\in C^\infty(M;\mathbb C)$ and $\eta\in\RX_x$ we have
$$
L((i_x^-)^*\alpha)(\eta+dh)
=L((i_x^-)^*(e^h\alpha))(\eta)e^{-h(x)}.
$$
The restriction $L((i_x^-)^*\alpha):\mathring\RX_x\to\mathbb C$ is
holomorphic.
If $\omega$ is Lyapunov for $X$ then $\RX_x+\omega\subseteq\mathring\RX_x$, 
and for all $\eta\in\RX_x$
\begin{equation}\label{E:cont:R}
\lim_{t\to0^+}L((i_x^-)^*\alpha)(\eta+t\omega)=L((i_x^-)^*\alpha)(\eta).
\end{equation}
Particularly, if $X$ satisfies \LL, then $\mathring\RX_x\subseteq\RX_x$ is dense,
and the function $L((i_x^-)^*\alpha):\RX_x\to\mathbb C$ is completely
determined by its restriction to $\mathring\RX_x$.
\end{proposition}

Be aware however, that without further assumptions the sets 
$\QX$, $\PX$ and $\RX$ might very well be empty.

\subsection{Morse complex and integration}\label{SS:Mc}

Let $\mathbb C^{\mathcal X}=\Maps(\mathcal X;\mathbb C)$ denote the vector space
generated by $\mathcal X$. Note that $\mathbb C^{\mathcal X}$ is $\mathbb Z$--graded
by $\mathbb C^{\mathcal X}=\bigoplus_q\mathbb C^{\mathcal X_q}$.
For $\eta\in \QX$ define a linear map
$$ 
\delta_\eta=\delta_\eta^{X,\mathcal O}:\mathbb C^{\mathcal X}\to\mathbb C^{\mathcal X},
\qquad
\delta_\eta(f)(x)
:=\sum_{y\in\mathcal X}L(\mathbb I_{x,y})(\eta)\cdot f(y)
$$
where $f\in\mathbb C^{\mathcal X}$ and $x\in\mathcal X$.
In section~\ref{SS:TUS} we will prove

\begin{proposition}\label{P:d2=0}
We have $\delta_\eta^2=0$, for all $\eta\in \QX$.
\end{proposition}

For a vector field $X$ which satisfies \LL\ and \MS, a choice
of orientations $\mathcal O$ and $\eta\in \QX$
we let $C^*_\eta(X;\mathbb C)=C^*_\eta(X,\mathcal O;\mathbb C)$ denote the complex with underlying
vector space $\mathbb C^{\mathcal X}$ and differential $\delta_\eta$.
Moreover, for $\eta\in\mathcal Z^1(M;\mathbb C)$ let $\Omega^*_\eta(M;\mathbb C)$ denote the deRham
complex with differential $d_\eta\alpha:=d\alpha+\eta\wedge\alpha$.
For $\eta\in\RX$ define a linear map 
$$
\Int_\eta=\Int_\eta^{X,\mathcal O}:\Omega^*(M;\mathbb C)\to\mathbb C^{\mathcal X},
\qquad
\Int_\eta(\alpha)(x):=L((i_x^-)^*\alpha)(\eta)
$$
where $\alpha\in\Omega^*(M;\mathbb C)$ and $x\in\mathcal X$.

The following two propositions will be proved in section~\ref{SS:TUS}.

\begin{proposition}\label{P:onto}
For $\eta\in \RX$ the linear map
$\Int_\eta:\Omega^*(M;\mathbb C)\to\mathbb C^{\mathcal X}$
is onto.
\end{proposition}

\begin{proposition}\label{P:inthom}
For $\eta\in \QX\cap \RX$ the integration is a 
homomorphism of complexes
\begin{equation}\label{E:inthom}
\Int_\eta:
\Omega^*_\eta(M;\mathbb C)\to C^*_\eta(X;\mathbb C).
\end{equation}
\end{proposition}

To make the gauge invariance more explicit,
suppose $h\in C^\infty(M;\mathbb C)$ and $\eta\in\QX\cap\RX$.
Then $\eta+dh\in\QX\cap\RX$, and we have a commutative
diagram of homomorphisms of complexes:
\begin{equation}\label{CD:gauge}
\xymatrix{
\Omega^*_\eta(M;\mathbb C)
\ar[rrr]^-{\Int_\eta}
&&& 
C^*_\eta(X;\mathbb C)
\\
\Omega^*_{\eta+dh}(M;\mathbb C)
\ar[u]^-{e^h}_-{\simeq} 
\ar[rrr]^-{\Int_{\eta+dh}}
&&& 
C^*_{\eta+dh}(X;\mathbb C)
\ar[u]_-{e^h}^-{\simeq}
}
\end{equation}

Let $\SX\subseteq\QX\cap\RX$ denote the subset of closed one forms
$\eta$ for which \eqref{E:inthom} does not induce an isomorphism in cohomology.
Note that $\SX$ is gauge invariant, i.e.\ $\SX+dh\subseteq\SX$
for $h\in C^\infty(M;\mathbb C)$.

Suppose $U$ is an open subset of a Fr\'echet space and let $S\subseteq U$ be
a subset. We say $S$ is an \emph{analytic subset} of $U$ if for every point
$z\in U$ there exists a neighborhood $V$ of $z$ and finitely many
holomorphic functions $f_1,\dotsc,f_N:V\to\mathbb C$ so that
$S\cap V=\{v\in V\mid f_1(v)=\cdots=f_N(v)=0\}$, see \cite{W72}.

\begin{theorem}\label{T:int}
Suppose $X$ satisfies \LL\ and \MS.
Then $\RX\subseteq\QX$. Moreover, $\SX\cap\mathring\RX$ is an
analytic subset of $\mathring \RX$.
If $\omega$ is a Lyapunov form for $X$ and $\eta\in\RX$, then 
there exists $t_0$ such that $\eta+t\omega\in\mathring\RX\setminus\Sigma$ for
all $t>t_0$. Particularly, the integration \eqref{E:inthom}
induces an isomorphism in cohomology for generic $\eta\in \RX$.
\end{theorem}

In general \eqref{E:inthom} will 
not induce an isomorphism in cohomology for all $\eta\in \RX$.
For example one can consider mapping cylinders and a nowhere vanishing $X$. 
In this case $\RX=\QX=\mathcal Z^1(M;\mathbb C)$, and
the complex $C^*_\eta(X;\mathbb C)$ is trivial.
However, the deRham cohomology is non-trivial for some $\eta$, e.g.\ $\eta=0$.

\subsection{Exponential growth}\label{SS:exp_grow}

In order to guaranty that $\RX$ is non-trivial we introduce

\begin{definition}[Exponential growth, \EG]\label{D:0} 
A vector field $X$ is said to have 
the \emph{exponential growth property at a rest point $x$} if for 
some (and then every) Riemannian metric $g$ on $M$ there exists 
$C\geq0$ so that $\Vol(B_x(r))\leq e^{Cr}$,
for all $r\geq0$. Here $B_x(r)\subseteq W^-_x$ denotes the  
ball of radius $r$ centered at $x\in W_x^-$ with respect
to the induced Riemannian metric $(i^-_x)^*g$ on $W_x^-$.
A vector field $X$ is said to have the \emph{exponential growth
property,} \EG\ for short, if it has the exponential growth property at all rest points.
\end{definition}

For rather trivial reasons every vector field with $\RX\neq\emptyset$ satisfies
\EG, see Proposition~\ref{P:RXEG}.
We are interested in the exponential growth property because of the
following converse statement which will be proved in section~\ref{SS:exp}.

\begin{proposition}\label{P:expR}
If $X$ satisfies \LL\ and \EG, then $\mathring\RX$ is non-empty. 
More precisely, if $\omega$ is a Lyapunov form for $X$ and $\eta\in\mathcal
Z^1(M;\mathbb C)$, then there exists $t_0\in\mathbb R$,
such that $\eta+t\omega\in\mathring\RX$ for all $t>t_0$.
\end{proposition}

\begin{remark}
Suppose $X$ satisfies \MS, \LL\ and \EG.
Let $\omega$ be a Lyapunov form for $X$, and let $x,y\in\mathcal X$.
In view of Proposition~\ref{P:expR} and Theorem~\ref{T:int} we have
$t\omega\in\QX_{x,y}$ for sufficiently large $t$. Hence
\begin{equation}\label{E:kkk}
\sum_{\gamma\in\mathcal P_{x,y}}\mathbb I_{x,y}(\gamma)e^{t\omega(\gamma)}
\end{equation}
converges absolutely for sufficiently large $t$. Particularly, there exists
$C\geq0$ such that $|\mathbb I_{x,y}(\gamma)|\leq e^{-C\omega(\gamma)}$, for
all $\gamma\in\mathcal P_{x,y}$. Since the sum \eqref{E:kkk} is
over homotopy classes, this is significantly stronger than what
was conjectured in \cite{N93} and proved in \cite{BH01} or \cite{P98}.
\end{remark}

Using a result of Pajitnov \cite{P98, P99, P03} we will prove the following
weak genericity result in section~\ref{SS:gen}.

\begin{theorem}\label{T:gen}  
Suppose $X$ satisfies \LL. Then, in every
$C^0$--neighborhood of $X$, there exists a vector field which
coincides with $X$ in a neighborhood of $\mathcal X$, and which satisfies
\LL, \MS, \NCT\ and \EG.
\end{theorem}

\begin{conjecture}\label{C:L}
If $X$ satisfies \LL, then in every $C^1$--neighborhood of $X$ there exists
a vector field which coincides with $X$ in a neighborhood of $\mathcal X$,
and which satisfies \LL, \MS, \NCT\ and \EG.
\end{conjecture}

For the sake of Theorem~\ref{T:tor} below we have to introduce the 
\emph{strong exponential growth property.}
Consider the bordism $W:=M\times[-1,1]$.
Set $\partial_\pm W:=M\times\{\pm1\}$.
Let $Y$ be a vector field on $W$.
Assume that there are vector fields $X_\pm$ on $M$ so that
$Y(z,s)=X_+(z)+(s-1)\partial/\partial s$
in a neighborhood of $\partial_+W$ and so that
$Y(z,s)=X_-(z)+(-s-1)\partial/\partial s$ in a neighborhood of $\partial_-W$.
Particularly, $Y$ is tangential to $\partial W$.
Moreover, assume that $ds(Y)<0$ on $M\times(-1,1)$.
Particularly, there are no zeros or closed trajectories of $Y$ contained
in the interior of $W$.
The properties \MS, \NCT, \LL\ and \EG\ make sense for these
kind of vector fields on $W$ too.

If $X$ satisfies \MS, \NCT\ and \LL, then it is easy to construct
a vector field $Y$ on $W$ as above satisfying \MS, \NCT\ and \LL\ 
such that $X_+=X$ and $X_-=-\grad_{g_0}f$ for a Riemannian metric $g_0$ 
on $M$ and a Morse function $f:M\to\mathbb R$, see Proposition~\ref{P:hl}
in appendix~\ref{app:B}. However, even if we assume that $X$
satisfies \EG, it is not clear that such a $Y$ can be chosen to have \EG.
We thus introduce the following, somewhat asymmetric,

\begin{definition}[Strong exponential growth, \SEG]\label{D:SEG}
A vector field $X$ on $M$ is said to have \emph{strong exponential growth}, \SEG\ for
short, if there exists a vector field $Y$ on $W=M\times[-1,1]$ as above
satisfying \MS, \NCT, \LL\ and \EG\
such that $X_+=X$ and $X_-=-\grad_{g_0}f$ for a Riemannian metric $g_0$ 
on $M$ and a Morse function $f:M\to\mathbb R$.
Note that \SEG\ implies \MS, \NCT, \LL\ and \EG.
\end{definition}

\begin{example}
A vector field without zeros satisfying \NCT\ and \LL\ satisfies \SEG. 
\end{example}

Using the same methods as for Theorem~\ref{T:gen} we will in
section~\ref{SS:gen} prove

\begin{theo_genp}
Suppose $X$ satisfies \LL. Then, in every $C^0$--neighborhood of $X$, there
exists a vector field which coincides with $X$ in a neighborhood of
$\mathcal X$ and satisfies \SEG.
\end{theo_genp}

\subsection{Torsion}\label{SS:tor}

Choose a Riemannian metric $g$ on $M$. Equip the space
$\Omega^*(M;\mathbb C)$ with a weakly non-degenerate bilinear form
$b(\alpha,\beta):=\int_M\alpha\wedge\star\beta$. For $\eta\in\Omega^1(M;\mathbb
C)$
let $d_\eta^t:\Omega^*(M;\mathbb C)\to\Omega^{*-1}(M;\mathbb C)$
denote the formal transpose of $d_\eta$ with respect to this bilinear form.
Explicitly, we have $d_\eta^t\alpha=d^*+i_{\sharp\eta}\alpha$,
where $\sharp\eta\in\Gamma(TM\otimes\mathbb C)$ is defined by
$g(\sharp\eta,\cdot)=\eta$. Consider the operator
$B_\eta=d_\eta\circ d_\eta^t+d_\eta^t\circ d_\eta$.
This is a zero order perturbation of the Laplace--Beltrami operator
and depends holomorphically on $\eta$.
Note that the adjoint of $B_\eta$ with respect to the standard Hermitian
structure on $\Omega^*(M;\mathbb C)$ coincides with
$B_{\bar\eta}$, where $\bar\eta$ denotes the complex conjugate of 
$\eta$.\footnote{This is called a `self adjoint holomorphic' family in \cite{K76}.}
Assume from now on that $\eta$ is closed. Then
$B_\eta$ commutes with $d_\eta$ and $d_\eta^t$.

For $\lambda\in\mathbb C$
let $E^*_\eta(\lambda)$ denote the generalized $\lambda$--eigen space of $B_\eta$.
Recall from elliptic theory that $E^*_\eta(\lambda)$ is finite dimensional
graded subspace $E_\eta^*(\lambda)\subseteq\Omega^*(M;\mathbb C)$.
The differentials $d_\eta$ and $d_\eta^t$ preserve $E^*_\eta(\lambda)$
since they commute with $B_\eta$. 
Note however that the restriction of $B_\eta-\lambda$ to
$E^*_\eta(\lambda)$ will in general only be nilpotent.
If $\lambda_1\neq\lambda_2$ then $E_\eta^*(\lambda_1)$
and $E_\eta^*(\lambda_2)$ are orthogonal with respect to $b$ since $B_\eta$
is symmetric with respect to $b$. It follows that 
$b$ restricts to a non-degenerate bilinear form on every $E_\eta^*(\lambda)$.
In section~\ref{SS:Tintii} we will prove

\begin{proposition}\label{P:hodge}
Let $\eta\in\mathcal Z^1(M;\mathbb C)$. Then
$E_\eta^*(\lambda)$ is acyclic for all $\lambda\neq0$, and
the inclusion $E^*_\eta(0)\to\Omega^*_\eta(M;\mathbb C)$ 
is a quasi isomorphism. 
\end{proposition}

If $\eta\in \RX\setminus\SX$ then, in view of Proposition~\ref{P:hodge}, 
the restriction of the integration
\begin{equation}\label{E:intE}
\Int_\eta|_{E^*_\eta(0)}:E^*_\eta(0)\to C^*_\eta(X;\mathbb C)
\end{equation}
is a quasi isomorphism. Recall that an endomorphism preserving
a non-degenerate bilinear form has determinant $\pm1$. Therefore $b$
determines an equivalence class of graded bases \cite{M66} in $E_\eta^*(0)$. Moreover, the
indicator functions provide a graded basis of $C^*_\eta(X;\mathbb C)$. 
Let $\pm T(\Int_\eta|_{E^*_\eta(0)})\in\mathbb C\setminus0$
denote the relative torsion of \eqref{E:intE} with respect to these bases,
see \cite{M66}. Moreover, define a complex valued Ray--Singer \cite{RS71}
kind of torsion
$$
(T^\an_{\eta,g})^2:=\prod_{q}\bigl(\det{}'B_\eta^q\bigr)^{(-1)^{q+1}q}
\in\mathbb C\setminus0
$$
where $\det{}'B_\eta^q$ denotes the zeta regularized product \cite{I01, BFKM96} of all
non-zero eigen values of $B_\eta^q:\Omega^q(M;\mathbb C)\to\Omega^q(M;\mathbb C)$,
computed with respect to the Agmon angle $\pi$.
In section~\ref{S:R} we will provide a regularization $R(\eta,X,g)$ of the
possibly divergent integral
$$
\int_{M\setminus\mathcal X}\eta\wedge X^*\Psi(g),
$$
where $\Psi(g)\in\Omega^{n-1}(TM\setminus M;\mathcal O_M)$ denotes the
global angular form.
Finally, set
\begin{equation}\label{E:TInt}
(T\Int_\eta)^2=
(T\Int^{X,\mathcal O}_{\eta,g})^2
:=(T(\Int_\eta|_{E_\eta^*(0)}))^2\cdot(T^\an_{\eta,g})^2
\cdot(e^{-R(\eta,X,g)})^2.
\end{equation}

In section~\ref{SS:PTint} we will show

\begin{proposition}\label{P:Tint}
The quantity \eqref{E:TInt} does not depend on $g$. It defines a function
\begin{equation}\label{E:Tintf}
(T\Int)^2:\RX\setminus\SX\to\mathbb C\setminus0
\end{equation}
which satisfies $(T\Int_{\bar\eta})^2=\overline{(T\Int_\eta)^2}$, and
which is gauge invariant, i.e.\ for $\eta\in\RX\setminus\Sigma$ and
$h\in C^\infty(M;\mathbb C)$ we have 
$$
(T\Int_{\eta+dh})^2=(T\Int_\eta)^2.
$$
The restriction $(T\Int)^2:\mathring\RX\setminus\Sigma\to\mathbb C\setminus0$
is holomorphic. If $\omega$ is Lyapunov for $X$ and $\eta\in\RX\setminus\Sigma$ then
for sufficiently small $t>0$ we have $\eta+t\omega\in\mathring\RX\setminus\Sigma$,
and 
\begin{equation}\label{E:cont:int}
\lim_{t\to0^+}(T\Int_{\eta+t\omega})^2=(T\Int_\eta)^2.
\end{equation}
\end{proposition}

\begin{remark}
In view of the gauge invariance $(T\Int)^2$ can be regarded as a partially
defined holomorphic function on the finite dimensional vector space
$H^1(M;\mathbb C)$.
\end{remark}

The rest of section~\ref{S:Ttor} is dedicated to the proof of

\begin{theorem}\label{T:tor}
Suppose $X$ satisfies \SEG.
Then $\mathring\PX$ is non-empty. More precisely, 
if $\omega$ is a Lyapunov form for $X$ and $\eta\in\mathcal Z^1(M;\mathbb C)$,
then there exists $t_0\in\mathbb R$ such that $\eta+t\omega\in\mathring\PX$
for all $t>t_0$. Moreover, for $\eta\in(\RX\setminus\SX)\cap\PX$
$$
(e^{L(h_*\mathbb P)(\eta)})^2=(T\Int_\eta)^2.
$$
Particularly, the zeta function $\eta\mapsto(e^{L(h_*\mathbb P)(\eta)})^2$ admits an
analytic continuation to $\mathring \RX$ with zeros and singularities contained in the proper
analytic subset $\mathring \RX\cap\SX$.
\end{theorem}

\begin{example}\label{Ex:maptor}
Let $f:N\to N$ be a diffeomorphism, and let $M$ denote the 
mapping cylinder obtained by glueing the boundaries of 
$N\times[0,1]$ with the help of $f$. Let $X=\partial/\partial t$, where
$t$ denotes the coordinate in $[0,1]$. 
Since it has no zeros at all $X$ satisfies \MS\ and \SEG.
Moreover, $X$ satisfies \NCT\ iff all fixed points of $f^k$ are non-degenerate
for all $k\in\mathbb N$. In this case we have
$$
e^{L(h_*\mathbb P)(zdt)}
=\exp\sum_{k=1}^\infty\sum_{x\in\Fix(f^k)}\frac{\ind_x(f^k)}k(e^z)^k
=\zeta_f(e^z)
$$
where $\zeta_f$ denotes the Lefschetz zeta function associated with $f$.
Theorem~\ref{T:tor} implies that for generic $z$ we have
$\pm T^\an_{zdt,g}=e^{zR(dt,X,g)}\zeta_f(e^z)$. This was already established by
Marcsik in his thesis \cite{M98}.
\end{example}

In the acyclic case it suffices to assume \EG.

\begin{theo_torp}
Suppose $X$ satisfies \LL, \MS, \NCT\ and \EG.
Assume that there exists $\eta_0\in\mathcal Z^1(M;\mathbb C)$
such that $H^*_{\eta_0}(M;\mathbb C)=0$.
Then $\mathring\PX$ is non-empty. More precisely, 
if $\omega$ is a Lyapunov form for $X$ and $\eta\in\mathcal Z^1(M;\mathbb C)$,
then there exists $t_0\in\mathbb R$ such that $\eta+t\omega\in\mathring\PX$
for all $t>t_0$. Moreover, for $\eta\in(\RX\setminus\SX)\cap\PX$
$$
(e^{L(h_*\mathbb P)(\eta)})^2=(T\Int_\eta)^2.
$$
Particularly, the zeta function $\eta\mapsto(e^{L(h_*\mathbb P)(\eta)})^2$ admits an
analytic continuation to $\mathring \RX$ with zeros and singularities contained in the proper
analytic subset $\mathring \RX\cap\SX$.
\end{theo_torp}

\begin{conjecture}
Theorem~\ref{T:tor}' remains true without the acyclicity assumption.
\end{conjecture}

\begin{remark}
Suppose $X$ satisfies \SEG. Let $\omega$ be 
a Lyapunov form for $X$. In view of Theorem~\ref{T:tor}
$$
\sum_{a\in H_1(M;\mathbb Z)/\Tor(H_1(M;\mathbb Z))}(h_*\mathbb P)(a)e^{t\omega(a)}
$$
converges absolutely for sufficiently large $t$. Particularly, there
exists $C\geq0$ such that $|(h_*\mathbb P)(a)|\leq e^{-C\omega(a)}$
for all $a\in H_1(M;\mathbb Z)/\Tor(H_1(M;\mathbb Z))$.
Note that for Pajitnov's class of vector fields the Laplace transform
$L(h_*\mathbb P)$ actually is a rational function \cite{P03}.
\end{remark}

\subsection{Interpretation via classical Dirichlet series}\label{SS:dir}

Restricting to affine lines $\eta+z\omega$ in $\mathcal Z^1(M;\mathbb C)$
we can interpret the above results in terms of classical Laplace transforms.

More precisely,
let $\eta\in\mathcal Z^1(M;\mathbb C)$, and suppose
$\omega$ is a Lyapunov form for $X$. If $X$ satisfies
\EG\ then there exists $\rho<\infty$ so that for all $x\in\mathcal X$ and
all $\alpha\in\Omega^*(M;\mathbb C)$ the Laplace transform
\begin{equation}\label{E:lapI}
\Int_{\eta+z\omega}(\alpha)(x)
=\int_0^\infty e^{-z\lambda}d\bigl((-h_x^\omega)_*(e^{h^\eta_x}(i_x^-)^*\alpha)\bigr)(\lambda)
\end{equation}
has abscissa of absolute convergence at most $\rho$, i.e.\ \eqref{E:lapI}
converges absolutely for all $\Re(z)>\rho$, see Proposition~\ref{P:expR}.
Here $(-h^\omega_x)_*(e^{h^\eta_x}(i_x^-)^*\alpha)$ denotes the push forward
of $e^{h^\eta_x}(i_x^-)^*\alpha$ considered as measure on $W_x^-$ via the
map $-h^\omega_x:W_x^-\to[0,\infty)$. The integral in \eqref{E:lapI} is
supposed to denote the Laplace transform of 
$(-h^\omega_x)_*(e^{h^\eta_x}(i_x^-)^*\alpha)$.

Assume in addition that $X$ satisfies \MS, and let $x,y\in\mathcal X$.
Consider the mapping $-\omega:\mathcal P_{x,y}\to\mathbb R$ 
and define a measure with discrete support,
$(-\omega)_*(\mathbb I_{x,y}e^\eta):[0,\infty)\to\mathbb C$ by
\begin{equation}\label{E:dirIc}
\bigl((-\omega)_*(\mathbb I_{x,y}e^\eta)\bigr)(\lambda)
:=\sum_{\{\gamma\in\mathcal P_{x,y}\mid -\omega(\gamma)=\lambda\}}
\mathbb I_{x,y}(\gamma)e^{\eta(\gamma)}.
\end{equation}
In view of Theorem~\ref{T:int} its Laplace transform,
i.e. the Dirichlet series
\begin{equation}\label{E:dirI}
L(\mathbb I_{x,y})(\eta+z\omega)=
\sum_{\lambda\in[0,\infty)}
e^{-z\lambda}\bigl((-\omega)_*(\mathbb I_{x,y}e^\eta)\bigr)(\lambda),
\end{equation}
has abscissa of absolute convergence at most $\rho$, i.e.\ \eqref{E:dirI}
converges absolutely for all $\Re(z)>\rho$.
Particularly, we see that from the germ at $+\infty$ of the holomorphic
function $z\mapsto\delta_{\eta+z\omega}$
one can recover, via inverse Laplace transform, a good amount of the 
counting functions $\mathbb I_{x,y}$, namely the numbers \eqref{E:dirIc}
for all $\lambda\in\mathbb R$ and all $x,y\in\mathcal X$.

Assume in addition that $X$ satisfies \NCT\ and \SEG. 
Consider the mapping $-\omega:[S^1,M]\to\mathbb R$
and define a measure with discrete support,
$(-\omega)_*(\mathbb Pe^\eta):[0,\infty)\to\mathbb C$ by
\begin{equation}\label{E:dirPc}
\bigl((-\omega)_*(\mathbb Pe^\eta)\bigr)(\lambda)
:=\sum_{\{\gamma\in[S^1,M]\mid -\omega(\gamma)=\lambda\}}
\mathbb P(\gamma)e^{\eta(\gamma)}.
\end{equation}
In view of
Theorem~\ref{T:tor} its Laplace transform,
i.e. the Dirichlet series
\begin{equation}\label{E:dirP}
L(h_*\mathbb P)(\eta+z\omega)=
\sum_{\lambda\in[0,\infty)}e^{-z\lambda}
\bigl((-\omega)_*(\mathbb Pe^\eta)\bigr)(\lambda),
\end{equation}
has finite abscissa of convergence, i.e.\ for sufficiently large $\Re(z)$
the series \eqref{E:dirP} converges absolutely.
Moreover, $z\mapsto e^{L(h_*\mathbb P)(\eta+z\omega)}$ admits an analytic
continuation with isolated singularities to
$\{z\in\mathbb C\mid\Re(z)>\rho\}$.
Particularly, we see that from the germ at $+\infty$ of the holomorphic
function $z\mapsto T\Int_{\eta+z\omega}$
one can recover, via inverse Laplace transform, a good amount of
the counting functions $\mathbb P$, namely the numbers 
\eqref{E:dirPc} for all $\lambda\in\mathbb R$.

\subsection{Relation with Witten--Helffer--Sj\"ostrand theory}

The above results provide some useful additions to
Witten--Helffer--Sj\"ostrand theory.
Recall that Witten--Helffer--Sj\"ostrand theory on a closed Riemannian manifold $(M,g)$ can be
extended from a Morse function $f:M\to\mathbb R$ to a
closed Morse one form $\omega\in\mathcal Z^1(M;\mathbb R)$, see \cite{BH01}. Precisely,
let $\eta\in\mathcal Z^1(M;\mathbb C)$ and consider the one parameter family of elliptic complexes
$\Omega^*_{\eta+t\omega}(M;\mathbb C):=(\Omega^*(M;\mathbb C),d_{\eta+t\omega})$,
equipped with the Hermitian scalar product induced by the Riemannian metric $g$.
Then, for sufficiently large $t$, we have a canonic orthogonal decomposition 
of cochain complexes
$$
\Omega^*_{\eta+t\omega}(M;\mathbb C)
=\Omega^*_{\eta+t\omega,\sm}(M;\mathbb C)\oplus
\Omega^*_{\eta+t\omega,\la}(M;\mathbb C).
$$

If $X$ is a smooth vector field with all the above properties including
exponential growth and having $\omega$ as a Lyapunov
closed one form then the restriction of the integration 
$$
\Int_{\eta+t\omega}:\Omega^*_{\eta+t\omega,\sm}(M;\mathbb C)
\to C^*_{\eta+t\omega}(X;\mathbb C)
$$
is an isomorphism for sufficiently large $t$.
In particular the canonical base of $\mathbb C^{\mathcal X}$ 
provides a canonical base 
$\{E_x(t)\}_{x\in\mathcal X}$ for the small complex 
$\Omega^*_{\eta+t\omega,\sm}(M;\mathbb C)$, 
and the differential $d_{\eta+t\omega}$,
when written in this base is a matrix whose components are the Laplace
transforms $L(\mathbb I_{x,y})(\eta+t\omega)$.
One can formulate this fact as: \emph{The counting of instantons is taken care
of by the small complex.}

The large complex $\Omega^*_{\eta+t\omega,\la}(M;\mathbb C)$ is acyclic and
has Ray--Singer torsion which in the case of a
Morse function $\omega=df$ is exactly $tR(\omega,X,g)+\log\Vol(t)$ where 
$\log\Vol(t):=\sum(-1)^q\log\Vol_q(t)$ and $\Vol_q(t)$
denotes the volume of the canonical base $\{E_x(t)\}_{x\in\mathcal X_q}$.
If $\omega$ is a non-exact form, the above
expression has an additional term $\Re(L(h_*\mathbb P)(\eta+t\omega))$.
One can formulate this fact as: \emph{The counting of closed trajectories is
taken care of by the large complex.}

\section{Exponential growth}\label{S:exp}

In section~\ref{SS:exp} we will reformulate the
exponential growth condition, see Proposition~\ref{P:expg2}, and
show that \EG\ implies $\mathring\RX\neq\emptyset$, i.e.\ prove Proposition~\ref{P:expR}.
In section~\ref{SS:virt} we will present a criterion which when satisfied
implies exponential growth, see Proposition~\ref{P:rhoxprop}. This 
criterion is satisfied by a class of vector fields 
introduced by Pajitnov. A theorem of Pajitnov tells that his class is
$C^0$--generic. Using this we will give a proof of
Theorem~\ref{T:gen} in section~\ref{SS:gen}.

\subsection{Exponential growth}\label{SS:exp}

Let $g$ be a Riemannian metric on $M$, and let $x\in\mathcal X$ be a zero of
$X$. Let $g_x:=(i_x^-)^*g$ denote the induced Riemannian metric on the
unstable manifold $W_x^-$. Let $r_x^g:=\dist_{g_x}(x,\cdot):W_x^-\to[0,\infty)$ 
denote the distance to $x$.
Let $B^g_x(s):=\{y\in W_x^-\mid r^g_x(y)\leq s\}$
denote the ball of radius $s$, and let
$\Vol_{g_x}(B^g_x(s))$ denote its volume.
Recall from Definition~\ref{D:0} that $X$ has the exponential growth
property at $x$ if there exists $C\geq0$ such that
$\Vol_{g_x}(B_x^g(s))\leq e^{Cs}$ for all $s\geq0$.
This does not depend on $g$ although $C$ does.

\begin{proposition}\label{P:expg2}
Let $X$ be a vector field and suppose $x\in\mathcal X$. Then $X$ has
exponential growth property at $x$ iff for one (and hence every)
Riemannian metric $g$ on $M$ there exists a constant $C\geq0$ such that
$e^{-Cr_x^g}\in L^1(W_x^-)$.
\end{proposition}

This proposition is an immediate consequence of the following two lemmas.

\begin{lemma}\label{L:exp:1}
Suppose there exists $C\geq0$ such that $\Vol_{g_x}(B_x^g(s))\leq e^{Cs}$ for all
$s\geq0$. Then
$e^{-(C+\epsilon)r_x^g}\in L^1(W_x^-)$ for every $\epsilon>0$.
\end{lemma}

\begin{proof}
Clearly
\begin{equation}\label{ee1}
\int_{W_x^-}e^{-(C+\epsilon)r_x^g}
=
\sum_{n=0}^\infty\int_{B_x^g(n+1)\setminus B_x^g(n)}e^{-(C+\epsilon)r^g_x}.
\end{equation}
On $B^g_x(n+1)\setminus B^g_x(n)$ we have $e^{-(C+\epsilon)r_x^g}\leq
e^{-(C+\epsilon)n}$ and thus
\begin{eqnarray*}
\int_{B^g_x(n+1)\setminus B^g_x(n)}e^{-(C+\epsilon)r_x^g}
&\leq&
\Vol_{g_x}(B^g_x(n+1))e^{-(C+\epsilon)n}   
\\&\leq&
e^{C(n+1)}e^{-(C+\epsilon)n}
=e^Ce^{-\epsilon n}
\end{eqnarray*}
So \eqref{ee1} implies
$$
\int_{W_x^-}e^{-(C+\epsilon)r_x^g}
\leq
e^C\sum_{n=0}^\infty e^{-\epsilon n} 
=e^C(1-e^{-\epsilon})^{-1}<\infty    
$$
and thus $e^{-(C+\epsilon)r_x^g}\in L^1(W_x^-)$.
\end{proof}

\begin{lemma}\label{L:exp:2}
Suppose we have $C\geq0$ such that $e^{-Cr_x^g}\in L^1(W_x^-)$.
Then there exists a constant $A>0$ such that
$\Vol_{g_x}(B_x^g(s))\leq Ae^{Cs}$ for all $s\geq0$.
\end{lemma}

\begin{proof}
We start with the following estimate for $N\in\mathbb N$:
\begin{align*}
\Vol_{g_x}(B^g_x(N+1))&e^{-C(N+1)}=
\\
&=\sum_{n=0}^N\Vol_{g_x}(B_x^g(n+1))e^{-C(n+1)}-\Vol_{g_x}(B_x^g(n))e^{-Cn}
\\ 
&\leq\sum_{n=0}^\infty\bigl(\Vol_{g_x}(B_x^g(n+1))-\Vol_{g_x}(B_x^g(n))\bigr)e^{-C(n+1)}
\\
&=\sum_{n=0}^\infty\Vol_{g_x}\bigl(B_x^g(n+1)\setminus B_x^g(n)\bigr)e^{-C(n+1)}
\\
&\leq\sum_{n=0}^\infty\int_{B_x^g(n+1)\setminus B_x^g(n)}e^{-Cr_x^g}
=\int_{W^-_x}e^{-Cr_x^g}
\end{align*}
Given $s\geq0$ we choose an integer $N$ with $N\leq s\leq N+1$. Then
$$
\Vol_{g_x}(B_x^g(s))e^{-Cs}\leq\Vol_{g_x}(B_x^g(N+1))e^{-CN}=e^C\Vol_{g_x}(B_x^g(N+1))e^{-C(N+1)},
$$
and thus $\Vol(B_x^g(s))e^{-Cs}
\leq e^C\int_{W_x^-}e^{-Cr_x^g}=:A<\infty$.
We conclude $\Vol_{g_x}(B_x^g(s))\leq Ae^{Cs}$ for all $s\geq0$.
\end{proof}

Let $\eta\in\mathcal Z^1(M;\mathbb C)$ be a closed one form.
Recall that $h^\eta_x:W_x^-\to\mathbb C$ denotes the unique smooth function
which satisfies $dh^\eta_x=(i_x^-)^*\eta$ and $h_x^\eta(x)=0$. Recall 
from section~\ref{SS:cf} that $\eta\in\RX_x$ if $e^{h_x^\eta}\in L^1(W_x^-)$.

\begin{proposition}\label{P:RXEG}
Let $X$ be a vector field and suppose $x\in\mathcal X$.
If $\RX_x\neq\emptyset$ then $X$ has exponential growth at $x$.
Particularly, if $\RX\neq\emptyset$ then $X$ satisfies \EG.
\end{proposition}

This proposition follows immediately from Proposition~\ref{P:expg2} and 
the following

\begin{lemma}\label{L:exp:4}
There exists a constant $C=C_{g,\eta}\geq0$ such that $|h_x^\eta|\leq Cr_x^g$.    
\end{lemma}

\begin{proof}
Suppose $y\in W_x^-$. For every path
$\gamma:[0,1]\to W_x^-$ with $\gamma(0)=x$ and $\gamma(1)=y$ we find
$$
|h_x^\eta(y)|
=
\Bigl|\int_0^1(dh_x^\eta)(\gamma'(t))dt\Bigr|
\leq
\|\eta\|\int_0^1|\gamma'(t)|dt
=
\|\eta\|\length(\gamma)
$$
where $\|\eta\|:=\sup_{z\in M}|\eta_z|_g$.
We conclude $|h_x^\eta(y)|\leq\|\eta\|r_x^g(y)$.
Hence we can take $C:=\|\eta\|$.
\end{proof}

Let us recall the following crucial estimate from \cite[Lemma~3]{BH01}.

\begin{lemma}\label{L:exp:3}
Suppose $\omega$ is a Lyapunov for $X$, and suppose $x\in\mathcal X$.
Then there exist $\epsilon=\epsilon_{g,\omega}>0$ and $C=C_{g,\omega}\geq0$ such that
$r_x^g\leq-Ch_x^\omega$ on $W_x^-\setminus B_x^g(\epsilon)$.
\end{lemma}

\begin{proof}[Proof of Proposition~\ref{P:expR}]
Suppose $x\in\mathcal X$.
In view of Lemma~\ref{L:exp:4} and Lemma~\ref{L:exp:3} there exists $C>0$
so that 
$\Re(h^\eta_x+th^\omega_x)\leq(C-t/C)r_x^g$
on $W_x^-\setminus B_x^g(\epsilon)$. Since $X$ has exponential growth at $x$
we have $|e^{h_x^{\eta+t\omega}}|=e^{\Re(h_x^\eta+th_x^\omega)}
\leq e^{(C-t/C)r_x^g}\in L^1(W_x^-)$, and hence $\eta+t\omega\in\RX_x$,
for sufficiently large $t$.
We conclude that $\eta+t\omega\in\mathring\RX=\bigcap_{x\in\mathcal X}\mathring\RX_x$ 
for sufficiently large $t$, see Proposition~\ref{P:basic:R}.
\end{proof}

\subsection{Virtual interactions}\label{SS:virt}

Suppose $N\subseteq M$ is an immersed submanifold of dimension $q$.
Let $\Gr_q(TM)$ denote the Grassmannified tangent bundle of $M$, i.e.\
the compact space of $q$--planes in $TM$.
The assignment $z\mapsto T_zN$ provides an immersion 
$N\subseteq\Gr_q(TM)$. We let $\overline{\Gr(N)}\subseteq\Gr_q(TM)$ denote 
the closure of its image. Moreover, for a zero $y\in\mathcal X$ we
let $\Gr_q(T_yW_y^-)\subseteq\Gr_q(TM)$ denote the Grassmannian of $q$--planes
in $T_yW_y^-$ considered as subset of $\Gr_q(TM)$.

\begin{definition}[Virtual interaction]\label{D:VI}
For a vector field $X$ and two zeros $x\in\mathcal X_q$ and $y\in\mathcal X$ we
define their \emph{virtual interaction} to be the 
compact set
$$
K_x(y):=\Gr_q(T_yW_y^-)\cap\overline{\Gr(W_x^-\setminus B)}
$$
where $B\subseteq W_x^-$ is a compact ball centered at $x$. Note that
$K_x(y)$ does not depend on the choice of $B$.
\end{definition}

Note that $K_x(y)$ is non-empty iff there exists a sequence $z_k\in W_x^-$
so that $\lim_{k\to\infty}z_k=y$ and so that $T_{z_k}W_x^-$
converges to a $q$--plane in $T_yM$ which is contained in $T_yW_y^-$.

Although we removed $B$ from $W_x^-$ the set
$K_x(x)$ might be non-empty. However, if we would not have removed $B$
the set $K_x(x)$ would never be empty for trivial reasons.
Because of dimensional reasons we have
$K_x(y)=\emptyset$ whenever $\ind(x)>\ind(y)$.
Moreover, it is easy to see that $K_x(y)=\emptyset$ whenever $\ind(y)=n$.

We are interested in virtual interactions because of the following

\begin{proposition}\label{P:rhoxprop}
Suppose $X$ satisfies \LL, let $x\in\mathcal X$, and assume
that the virtual interactions $K_x(y)=\emptyset$ for all $y\in\mathcal X$. 
Then $X$ has exponential growth at $x$.
\end{proposition}

To prove Proposition ~\ref{P:rhoxprop} we will need the following

\begin{lemma}\label{little_lemm}
Let $(V,g)$ be an Euclidean vector space and $V=V^+\oplus V^-$ an orthogonal
decomposition. For $\kappa\geq0$ consider the endomorphism
$A_\kappa:=\kappa\id\oplus-\id\in\End(V)$ and the function
$$
\delta^{A_\kappa}:\Gr_q(V)\to\mathbb R,\quad
\delta^{A_\kappa}(W):=\tr_{g|_W}(p^\perp_W\circ A_\kappa\circ i_W),
 $$
where $i_W:W\to V$ denotes the inclusion and $p^\perp_W:V\to W$ the
orthogonal projection. Suppose we have a compact subset $K\subseteq\Gr_q(V)$
for which $\Gr_q(V^+)\cap K=\emptyset$. Then there exists $\kappa>0$ and
$\epsilon>0$ with $\delta^{A_\kappa}\leq-\epsilon$ on $K$.
\end{lemma}

\begin{proof}
Consider the case $\kappa=0$. Let $W\in\Gr_q(V)$ and choose a
$g|_W$ orthonormal base $e_i=(e^+_i,e^-_i)\in V^+\oplus V^-$,
$1\leq i\leq q$, of $W$. Then
$$
\delta^{A_0}(W)=\sum_{i=1}^qg(e_i,A_0e_i)
=-\sum_{i=1}^qg(e_i^-,e_i^-).
$$
So we see that $\delta^{A_0}\leq 0$ and
$\delta^{A_0}(W)=0$ iff $W\in\Gr_q(V^+)$. Thus $\delta^{A_0}|_K<0$.
Since $\delta^{A_\kappa}$ depends continuously on $\kappa$ and since $K$ is
compact we certainly find $\kappa>0$ and $\epsilon>0$ so that
$\delta^{A_\kappa}|_K\leq-\epsilon$.
\end{proof}

\begin{proof}[Proof of Proposition~\ref{P:rhoxprop}]
Let $S\subseteq W_x^-$ denote a small sphere centered at $x$.
Let $\tilde X:=(i_x^-)^*X$ denote the restriction of $X$ to $W_x^-$ and
let $\Phi_t$ denote the flow of $\tilde X$ at time $t$. Then
$$
\varphi:S\times[0,\infty)\to W^-_x,\quad
\varphi(x,t)=\varphi_t(x)=\Phi_t(x)
$$
parameterizes $W_x^-$ with a small neighborhood of $x$ removed.

Let $\kappa>0$. For every $y\in\mathcal X$ choose a chart
$u_y:U_y\to\mathbb R^n$ centered at $y$ so that
$$
X|_{U_y}=\kappa
\sum_{i\leq\ind(y)}u_y^i\frac\partial{\partial u_y^i}
-\sum_{i>\ind(y)}u_y^i\frac\partial{\partial u_y^i}.
$$
Let $g$ be a Riemannian metric on $M$ which restricts to
$\sum_idu_y^i\otimes du_y^i$ on $U_y$ and set $g_x:=(i_x^-)^*g$. Then
$$
\nabla X|_{U_y}=\kappa
\sum_{i\leq\ind(y)}du_y^i\otimes\frac\partial{\partial u_y^i}
-\sum_{i>\ind(y)}du_y^i\otimes\frac\partial{\partial u_y^i}.
$$
In view of our assumption $K_x(y)=\emptyset$ for all $y\in\mathcal X$
Lemma~\ref{little_lemm} permits us to choose $\kappa>0$ and $\epsilon>0$ so
that after possibly shrinking $U_y$ we have
\begin{equation}\label{divestia}
\diver_{g_x}(\tilde X)
=\tr_{g_x}(\nabla\tilde X)
\leq-\epsilon<0
\quad\text{on}\quad
\varphi(S\times[0,\infty))\cap(i_x^-)^{-1}
\Bigl(\bigcup_{y\in\mathcal X}U_y\Bigr).
\end{equation}
Let $\omega$ be a Lyapunov form for $X$. Since $\omega(X)<0$
on $M\setminus\mathcal X$, we can choose $\tau>0$ so that
\begin{equation}\label{divestib}
\tau\omega(X)+\ind(x)||\nabla X||_g\leq-\epsilon<0
\quad\text{on}\quad
M\setminus\bigcup_{y\in\mathcal X}U_y.
\end{equation}
Using $\tau\tilde X\cdot h^\omega_x\leq 0$ and
$$
\diver_{g_x}(\tilde X)
=\tr_{g_x}(\nabla\tilde X)
\leq\ind(x)||\nabla\tilde X||_{g_x}
\leq\ind(x)||\nabla X||_g
$$
\eqref{divestia} and \eqref{divestib} yield
\begin{equation}\label{divesti}
\tau \tilde X\cdot h^\omega_x+\diver_{g_x}(\tilde X)\leq-\epsilon<0
\quad\text{on}\quad
\varphi(S\times[0,\infty)).
\end{equation}
Choose an orientation of $W_x^-$ and let $\mu$ denote the volume form
on $W_x^-$ induced by $g_x$. Consider the function
$$
\psi:[0,\infty)\to\mathbb R,\quad
\psi(t):=\int_{\varphi(S\times[0,t])}e^{\tau h^\omega_x}\mu\geq 0.
$$
For its first derivative we find
$$
\psi'(t)=\int_{\varphi_t(S)}e^{\tau h^\omega_x}i_{\tilde X}\mu>0
$$
and for the second derivative, using \eqref{divesti},
\begin{eqnarray*}
\psi''(t)
&=&
\int_{\varphi_t(S)}
\bigl(\tau\tilde X\cdot h^\omega_x+\diver_{g_x}(\tilde X)\bigr)
e^{\tau h^\omega_x}i_{\tilde X}\mu
\\
&\leq&
-\epsilon\int_{\varphi_t(S)}e^{\tau h^\omega_x}i_{\tilde X}\mu
=-\epsilon\psi'(t).
\end{eqnarray*}
So $(\ln\circ\psi')'(t)\leq-\epsilon$ hence
$\psi'(t)\leq\psi'(0)e^{-\epsilon t}$
and integrating again we find
$$
\psi(t)
\leq\psi(0)+\psi'(0)(1-e^{-\epsilon t})/\epsilon
\leq\psi'(0)/\epsilon.
$$
So we have $e^{\tau h^\omega_x}\in L^1\bigl(\varphi(S\times[0,\infty)\bigr)$
and hence $e^{\tau h^\omega_x}\in L^1(W_x^-)$ too. We conclude
$\tau\omega\in\RX_x$. From Proposition~\ref{P:RXEG}
we see that $X$ has exponential growth at $x$.
\end{proof}


\subsection{Proof of Theorem~\ref{T:gen}}\label{SS:gen}

Let $X$ be a vector field satisfying \LL. Using Proposition~\ref{P:cone} we 
find a Lyapunov form $\omega$ for $X$ with integral cohomology class.
Hence there exists a smooth function $\theta:M\to S^1$ so that 
$\omega=d\theta$ is Lyapunov for $X$.

Choose a regular value $s_0\in S^1$ of $\theta$.
Set $V:=\theta^{-1}(s_0)$ and let $W$ denote the bordism obtained by cutting
$M$ along $V$, i.e.\ $\partial_\pm W=V$. This construction
provides a diffeomorphism $\Phi:\partial_-W\to\partial_+W$. Such a pair
$(W,\Phi)$ is called a cyclic bordism in \cite{P99}. When referring to Pajitnov's work
below we will make precise references to \cite{P99} but see also
\cite{P98} and \cite{P03}.

We continue to denote by $X$ the vector field on
$W$ induced from $X$, and by $\theta:W\to[0,1]$ the map induced from $\theta$.
We are exactly in the situation of Pajitnov: $-X$ is a
$\theta$--gradient in the sense of \cite[Definition~2.3]{P99}.
In view of \cite[Theorem~4.8]{P99}
we find, arbitrarily $C^0$--close to $X$, a smooth vector field $Y$ on $W$
which coincides with $X$ in a neighborhood of $\mathcal X\cup\partial W$, 
and so that $-Y$ is a $\theta$--gradient
satisfying condition $(\mathfrak C\mathcal Y)$ from \cite[Definition~4.7]{P99}.
For the reader's convenience we will below review Pajitnov's condition
$(\mathfrak C\mathcal Y)$ in more details.

Since $X$ and $Y$ coincide in a neighborhood of $\partial W$, $Y$
defines a vector field on $M$ which we denote by $Y$ too.
Clearly, $\omega=d\theta$ is Lyapunov for $Y$.
Using the $C^0$--openness statement in \cite[Theorem~4.8]{P99} 
and Proposition~\ref{P:KS} we may, by performing a $C^1$--small 
perturbation of $Y$, assume that $Y$ in addition satisfies \MS\ and \NCT.
Obviously, condition $(\mathfrak C\mathcal Y)$ implies that
$K^Y_x(y)=\emptyset$ whenever $\ind(x)\leq\ind(y)$, see below.
For trivial reasons we have $K^Y_x(y)=\emptyset$
whenever $\ind(x)>\ind(y)$.
It now follows from Proposition~\ref{P:rhoxprop} that $Y$ satisfies \EG\ too. 
This completes the proof of Theorem~\ref{T:gen}.

We will now turn to Pajitnov's condition $(\mathfrak C\mathcal Y)$, see
\cite[Definition~4.7]{P99}. Recall first that a smooth vector field $-X$ on a
closed manifold $M$ which satisfies \MS\ and is an $f$--gradient in the sense
of \cite[Definition~2.3]{P99} for some Morse function $f$, provides a partition of the 
manifold in cells, the unstable sets of the rest points of $-X$.
We will refer to such a partition as a \emph{generalized triangulation.}
The union of the unstable sets of $-X$ of rest points of index at most $k$
represents the $k$--skeleton and will be denoted \cite[section~2.1.4]{P99} by
$$
D(\ind\leq k,-X).
$$
>From this perspective the dual triangulation is associated to the
vector field $X$ which has the same properties with respect to $-f$.

Given a Riemannian metric $g$ on $M$ we will also write 
$$
B_\delta (\ind\leq k,-X)
\quad\text{resp.}\quad
D_\delta (\ind\leq k,-X)
$$
for the open resp.\ closed $\delta$--thickening of
$D(\ind\leq k,-X)$. They are the sets of points which lie on
trajectories of $-X$ which depart from the open resp.\ closed ball of 
radius $\delta$ centered at the rest points of Morse index at most $k$. 
It is not hard to see \cite[Proposition~2.30]{P99} that when $\delta\to0$ the
sets $B_\delta (\ind\leq k,-X)$ resp.\
$D_\delta(\ind\leq k,-X)$ provide a fundamental system of open resp.\
closed neighborhoods of $D(\ind\leq k,-X)$.
We also write
$$
C_\delta(\ind\leq k,-X):= M\setminus B_\delta(\ind\leq n-k-1,X).
$$
Note that for sufficiently small $\delta>0$ 
$$
B_\delta(\ind\leq k,-X)\subseteq C_\delta(\ind\leq k,-X).
$$

These definitions and notations can be also used in the case of a 
bordism, see \cite{P99} and \cite{M66}.
Denote by $U_\pm\subseteq\partial_\pm W$ the set of points $y\in\partial_\pm W$ 
so that the trajectory of the vector field $-X$ trough $y$ arrives resp.\
departs from $\partial W_\mp$ at some positive resp.\ negative time $t$.  
They are open sets. Following Pajitnov's notation we denote by 
$(-X)^\rightsquigarrow:U_+\to U_-$ resp.\
$X^\rightsquigarrow:U_-\to U_+$ the obvious diffeomorphisms induced by 
the flow of $X$ which are inverse one to the other.
If $A\subseteq\partial_\pm W$ we write for simplicity $(\mp X)^\rightsquigarrow(A)$ 
instead of $(\mp X)^\rightsquigarrow(A\cap U_\pm)$.

\begin{definition}[Property $(\mathfrak C\mathcal Y)$, see
{\cite[Definition~4.7]{P99}}]
A gradient like vector field $-X$ on a cyclic bordism $(W,\Phi)$ 
satisfies $(\mathfrak C\mathcal Y)$ if
there exist generalized triangulations $X_\pm$ on $\partial_\pm W$
and sufficiently small $\delta>0$ so that the following hold:
$$
\Phi(X_-)=X_+
$$
\begin{multline}\label{E:B2}
\tag{B$+$}
X^\rightsquigarrow\bigl(C_\delta(\ind\leq k,X_-)\bigr)
\cup\bigl(D_\delta(\ind\leq k+1,X)\cap\partial_+W\bigr)
\\
\subseteq B_\delta(\ind\leq k,X_+)
\end{multline}
\begin{multline}\label{E:B1}
\tag{B$-$}
(-X)^\rightsquigarrow\bigl(C_\delta(\ind\leq k,-X_+)\bigr)
\cup\bigl(D_\delta(\ind\leq k+1,-X)\cap\partial_-W\bigr)
\\
\subseteq B_\delta(\ind\leq k,-X_-)
\end{multline}
\end{definition}

If the vector field $Y$ on $(W,\Phi)$ constructed by the cutting off construction
satisfies $(\mathfrak C\mathcal Y)$ then, when regarded on $M$, 
it has the following property: Every zero $y$ admits a neighborhood
which does not intersect the unstable set of a zero $x$ with $\ind(y)\geq\ind(x)$.
Hence the virtual interaction $K^Y_x(y)$ is empty.
This is exactly what we used in the derivation of Theorem~\ref{T:gen} above.

Using Proposition~\ref{P:hl} in appendix~\ref{app:B} and Proposition~\ref{P:cone}
it is a routine task to extend the considerations 
above to the manifold $M\times[-1,1]$
and prove Theorem~\ref{T:gen}' along the same lines.

\section{The regularization $R(\eta,X,g)$}\label{S:R}

In this section we discuss the numerical invariant $R(\eta,X,g)$ associated
with a vector field $X$, a closed one form $\eta\in\mathcal Z^1(M;\mathbb C)$ and 
a Riemannian metric $g$. The invariant is defined by a possibly divergent 
but regularizable integral. It is implicit in the work of Bismut--Zhang \cite{BZ92}. 
More on this invariant is contained in \cite{BH03}.

Throughout this section we assume that $M$ is a closed manifold of dimension
$n$, and $X$ is a smooth vector field with zero set $\mathcal X$. We 
assume that the zeros are non-degenerate but not necessarily of the form
\eqref{E:4}. It is not difficult to generalize the regularization
to vector fields with isolated singularities, see \cite{BH03}.

\subsection{Euler, Chern--Simons, and the global angular form}\label{SS:MQECS}

Let $\pi:TM\to M$ denote the tangent bundle, and let $\mathcal O_M$ denote
the orientation bundle, a flat real line bundle over $M$. 
For a Riemannian metric $g$ let
$$
\ec(g)\in\Omega^n(M;\mathcal O_M)
$$
denote its \emph{Euler form}. For two Riemannian metrics $g_1$ and $g_2$ let
$$
\cs(g_1,g_2)\in
\Omega^{n-1}(M;\mathcal O_M)/d(\Omega^{n-2}(M;\mathcal O_M))
$$ 
denote their \emph{Chern--Simons class.} The definition of both quantities is
implicit in the formulae \eqref{dpsie} and \eqref{pg1pg2cs} below.
They have the following properties which follow immediately from
\eqref{dpsie} and \eqref{pg1pg2cs} below.
\begin{eqnarray}
d\cs(g_1,g_2) &=& \ec(g_2)-\ec(g_1)  
\label{E:csg:i}
\\
\cs(g_2,g_1) &=& -\cs(g_1,g_2)   
\label{E:csg:ii}
\\
\cs(g_1,g_3) &=& \cs(g_1,g_2) + \cs(g_2,g_3)
\label{E:csg:iii}
\end{eqnarray}

Let $\xi$ denote the Euler vector field on $TM$ which
assigns to a point $x\in TM$ the vertical vector $-x\in T_xTM$.
A Riemannian metric $g$ determines the Levi--Civita connection in the bundle
$\pi:TM\to M$. There is a canonic $\vvol(g)\in\Omega^n(TM;\pi^*\mathcal O_M)$,
which vanishes when contracted with horizontal vectors and which assigns
to an $n$--tuple of vertical vectors \emph{their volume times their
orientation.} The global angular form, see for instance \cite{BT82},
is the differential form
$$
\Psi(g)
:=\frac{\Gamma(n/2)}{(2\pi)^{n/2}|\xi|^n}
i_\xi\vvol(g)\in\Omega^{n-1}(TM\setminus M;\pi^*\mathcal O_M).
$$
This form was also considered by Mathai and Quillen \cite{MQ86}, and
was referred to as the Mathai--Quillen form in \cite{BZ92}.
Note that $\Psi(g)$ is the pull back of a form on $(TM\setminus M)/\R_+$.
Moreover, we have the equalities:
\begin{eqnarray}
d\Psi(g) &=& \pi^*\ec(g).
\label{dpsie}
\\
\Psi(g_2)-\Psi(g_1) &=& \pi^*\cs(g_1,g_2) \mod\pi^*d\Omega^{n-2}(M;\mathcal O_M)
\label{pg1pg2cs}
\end{eqnarray}
Further, if $x\in\mathcal X$ then
\begin{equation}\label{psiind}
\lim_{\epsilon\to0}\int_{\partial(M\setminus B_x(\epsilon))}
X^*\Psi(g)=\IND(x),
\end{equation}
where $\IND(x)$ denotes the \emph{Hopf index} of $X$ at $x$, and $B_x(\epsilon)$
denotes the ball of radius $\epsilon$ centered at $x$.

\subsection{Euler and Chern--Simons class for vector fields}\label{SS:ECSX}

Let $C_k(M;\mathbb Z)$ denote the complex of smooth singular chains in $M$.
Define a singular zero chain
$$
\ec(X):=\sum_{x\in\mathcal X}\IND(x)x\in C_0(M;\mathbb Z).
$$
For two vector fields $X_1$ and $X_2$ we are going to define 
\begin{equation}\label{E:csx}
\cs(X_1,X_2)\in C_1(M;\mathbb Z)/\partial C_2(M;\mathbb Z)
\end{equation}
with the following properties analogous to 
\eqref{E:csg:i}--\eqref{E:csg:iii}.
\begin{eqnarray}
\partial\cs(X_1,X_2) &=& \ec(X_2)-\ec(X_1)
\label{E:csx:i}
\\
\cs(X_2,X_1) &=& -\cs(X_1,X_2)
\label{E:csx:ii}
\\
\cs(X_1,X_3) &=& \cs(X_1,X_2) + \cs(X_2,X_3)
\label{E:csx:iii}    
\end{eqnarray}

It is constructed as follows.
Consider the vector bundle $p^*TM\to I\times M$, where $I:=[1,2]$ and
$p:I\times M\to M$ denotes the natural projection. Choose a section
$\mathbb X$ of $p^*TM$ which is transversal to the zero section and which
restricts to $X_i$ on $\{i\}\times M$, $i=1,2$. The zero set of $\mathbb X$
is a canonically oriented one dimensional submanifold with boundary. Its
fundamental class, when pushed forward via $p$, gives rise to
$c(\mathbb X)\in C_1(M;\mathbb Z)/\partial C_2(M;\mathbb Z)$.
Clearly $\partial c(\mathbb X)=\ec(X_2)-\ec(X_1)$.

Suppose $\mathbb X_1$ and $\mathbb X_2$ are two non-degenerate homotopies
from $X_1$ to $X_2$. Then
$c(\mathbb X_1)=c(\mathbb X_2)\in C_1(M;\mathbb Z)/\partial C_2(M;\mathbb Z)$.
Indeed, consider the vector bundle
$q^*TM\to I\times I\times M$, where $q:I\times I\times M\to M$ denotes the
natural projection. Choose a section of $q^*TM$ which is transversal to the
zero section, restricts to $\mathbb X_i$ on $\{i\}\times I\times M$,
$i=1,2$, and which restricts to $X_i$ on $\{s\}\times\{i\}\times M$ for all
$s\in I$
and $i=1,2$. The zero set of such a section then gives rise to $\sigma$
satisfying $c(\mathbb X_2)-c(\mathbb X_1)=\partial\sigma$.
Hence we may define $\cs(X_1,X_2):=c(\mathbb X)$.

\subsection{The regularization}\label{SS:reg}

Let $g$ be a Riemannian metric, and let
$\eta\in\mathcal Z^1(M;\mathbb C)$.
Choose a smooth function $f:M\to\mathbb C$ so that 
$\eta':=\eta-df$ vanishes on a neighborhood of 
$\mathcal X$. Then the following expression is well defined:
\begin{equation}\label{R:6}
R(\eta,X,g;f)
:=\int_{M\setminus\mathcal X}\eta'\wedge X^*\Psi(g)
-\int_Mf\ec(g)
+\sum_{x\in\mathcal X}\IND(x)f(x)
\end{equation}

\begin{lemma}\label{L:Rind}
The quantity $R(\eta,X,g;f)$ is independent of $f$.
\end{lemma}

\begin{proof}
Suppose $f_1$ and $f_2$ are two functions such that 
$\eta'_i:=\eta-df_i$ vanishes in a neighborhood of 
$\mathcal X$, $i=1,2$. Then $f_2-f_1$ is locally 
constant near $\mathcal X$. Using \eqref{psiind} and Stokes'
theorem we therefore get
$$
\int_{M\setminus\mathcal X}d\bigl((f_2-f_1)X^*\Psi(g)\bigr)
=\sum_{x\in\mathcal X}\bigl(f_2(x)-f_1(x)\bigr)\IND(x).
$$
Together with \eqref{dpsie} this immediately implies
$R(\eta,X,g;f_1)=R(\eta,X,g;f_2)$.
\end{proof}

\begin{definition}\label{D:ReXg}
For a vector field $X$ with non-degenerate zeros, a Riemannian metric $g$
and a closed one form $\eta\in\Omega^1(M;\mathbb C)$ we define
$R(\eta,X,g)$ by \eqref{R:6}. In view of 
Lemma~\ref{L:Rind} this does not depend on the choice of $f$.
We think of $R(\eta,X,g)$ as regularization of the possibly divergent
integral $\int_{M\setminus\mathcal X}\eta\wedge X^*\Psi(g)$.
\end{definition}

\begin{proposition}\label{P:Rdh}
For a smooth function $h:M\to\mathbb C$ we have
\begin{equation}\label{E:Rdh}
R(\eta+dh,X,g)-R(\eta,X,g)=-\int_Mh\ec(g)+\sum_{x\in\mathcal X}\IND(x)h(x).
\end{equation}
\end{proposition}

\begin{proof}
This is trivial, $h$ can be absorbed in the choice of $f$.
\end{proof}

\begin{proposition}\label{P:rg1g2}
For two Riemannian metrics $g_1$ and $g_2$ we have
\begin{equation}\label{E:rg1g2}
R(\eta,X,g_2)-R(\eta,X,g_1)
=\int_M\eta\wedge\cs(g_1,g_2).
\end{equation}
\end{proposition}

\begin{proof}
This follows easily from \eqref{pg1pg2cs}, Stokes' theorem and 
\eqref{E:csg:i}.
\end{proof}

\begin{proposition}\label{P:38}
For two vector fields $X_1$ and $X_2$ we have
\begin{equation}\label{EE:38}
R(\eta,X_2,g)-R(\eta,X_1,g)=\eta(\cs(X_1,X_2)).
\end{equation}
\end{proposition}

\begin{proof}
In view of \eqref{E:Rdh} and \eqref{E:csx:i}
we may w.l.o.g.\ assume that $\eta$ vanishes on a neighborhood of 
$\mathcal X_1\cup\mathcal X_2$.
Choose a non-degenerate homotopy $\mathbb X$ from $X_1$ to $X_2$.
Perturbing the homotopy, cutting it into several pieces and using
\eqref{E:csx:iii}
we may further assume that the zero set $\mathbb X^{-1}(0)\subseteq I\times M$
is actually contained in a simply connected $I\times V$. Again, we may
assume that $\eta$ vanishes on $V$. Then the right hand side of
\eqref{EE:38} obviously vanishes. Moreover, in
this situation Stokes' theorem implies
\begin{align*}
R(\eta,X_2,g)-R(\eta,X_1,g)
&=\int_{M\setminus V}\eta\wedge X_2^*\Psi(g)
-\int_{M\setminus V}\eta\wedge X_1^*\Psi(g)
\\&=
\int_{I\times(M\setminus V)}d(p^*\eta\wedge\mathbb X^*\tilde p^*\Psi(g))
\\&=
-\int_{I\times(M\setminus V)}p^*(\eta\wedge\ec(g))=0.
\end{align*}
Here $p:I\times M\to M$ denotes the natural projection, and
$\tilde p:p^*TM\to TM$ denotes the natural vector bundle homomorphism
over $p$. For the last calculation note that $d\mathbb X^*\tilde
p^*\Psi(g)=p^*\ec(g)$
in view of \eqref{dpsie}, and that $\eta\wedge\ec(g)=0$
because of dimensional reasons.
\end{proof}

\section{Completion of trajectory spaces and unstable manifolds}

If a vector field satisfies \MS\ and \LL, then the space of trajectories as well
as the unstable manifolds can be completed to manifolds with corners.
In section~\ref{SS:TUS} we recall these results, see Theorem~\ref{T:6}
below, and use them to prove Propositions~\ref{P:d2=0} and \ref{P:inthom}.
The rest of this section is dedicated to the proof of Theorem~\ref{T:int}.

\subsection{The completion}\label{SS:TUS}

Let $X$ be vector field on the closed manifold $M$ and suppose that 
$X$ satisfies \MS. Let $\pi:\tilde M\to M$ denote the universal covering.
Denote by $\tilde X$ the vector field $\tilde X:=\pi^*X$
and set $\tilde{\mathcal X}=\pi^{-1}(\mathcal X)$.

Given $\tilde x\in\tilde{\mathcal X}$
let $i^\pm_{\tilde x}:W^\pm_{\tilde x}\to\tilde M$ denote the
one-to-one immersions whose images define the stable and
unstable sets of $\tilde x$ with respect to the vector field
$\tilde X$. 
For any $\tilde x$ with $\pi(\tilde x)=x$ one can
canonically identify $W^\pm_{\tilde x}$ to $W^\pm_x$
so that $\pi\circ i^\pm_{\tilde x}=i^\pm_x$.
Define $\mathcal M(\tilde x,\tilde y):=
W^-_{\tilde x}\cap W^+_{\tilde y}$ if $\tilde x\neq\tilde y$,
and set $\mathcal M(\tilde x,\tilde x):=\emptyset$.
As the maps $i^-_{\tilde x}$ and $ i^+_{\tilde y}$ are
transversal, $\mathcal M(\tilde x,\tilde y)$
is a submanifold of $\tilde M$ of dimension 
$\ind(\tilde x)-\ind(\tilde y)$. It is
equipped with a free $\mathbb R$--action
defined by the flow generated by $\tilde X$.
Denote the  quotient
$\mathcal M(\tilde x,\tilde y)/\mathbb R$ by
$\mathcal T(\tilde x,\tilde y)$.
The quotient $\mathcal T(\tilde x,\tilde y)$ is a smooth
manifold of dimension
$\ind(\tilde x)-\ind(\tilde y)-1$, possibly empty.
If $\ind(\tilde x)\leq\ind(\tilde y)$,
in view the transversality required
by the hypothesis \MS, the manifolds
$\mathcal M(\tilde x,\tilde y)$ and
$\mathcal T(\tilde x,\tilde y)$ are empty.

An \emph{unparameterized broken trajectory} from
$\tilde x\in\tilde{\mathcal X}$
to $\tilde y\in\tilde{\mathcal X}$, is an element of the set
$\hat{\mathcal T}(\tilde x,\tilde y)
:=\bigcup_{k\geq0}\hat{\mathcal T}(\tilde x,\tilde y)_k$, where
\begin{equation}\label{Bk_defi}
\hat{\mathcal T}(\tilde x,\tilde y)_k
:=\bigcup
\mathcal T(\tilde y_0,\tilde y_1)
\times\cdots\times\mathcal T(\tilde y_k,\tilde y_{k+1})
\end{equation}
and the union is over all (tuples of) critical points
$\tilde y_i\in\tilde{\mathcal X}$ with $\tilde y_0=\tilde x$
and $\tilde y_{k+1}=\tilde y$.

For $\tilde x\in \tilde{\mathcal X}$ introduce the
\emph{completed unstable set} $\hat W^-_{\tilde x}:=
\bigcup_{k\geq0}(\hat W^-_{\tilde x})_k$, where
\begin{equation}\label{Wk_defi}
(\hat W^-_{\tilde x})_k
:=\bigcup
\mathcal T(\tilde y_0,\tilde y_1)\times\cdots\times
\mathcal T(\tilde y_{k-1},\tilde y_k)\times W^-_{\tilde y_k}
\end{equation}
and the union is over all (tuples of) critical points
$\tilde y_i\in\tilde{\mathcal X}$ with $\tilde y_0=\tilde x$.

Let $\hat i_{\tilde x}^-:\hat W^-_{\tilde x}\to\tilde M$ denote
the map whose restriction to  
$\mathcal T(\tilde y_0,\tilde y_1)
\times\cdots\times\mathcal T(\tilde y_{k-1},\tilde y_k)
\times W^-_{\tilde y_k}$
is the composition of the projection on $W^-_{\tilde y_k}$ with
$i_{\tilde y_k}^-$.

Recall that an $n$--dimensional \emph{manifold with corners}
$P$, is a paracompact Hausdorff space equipped with a
maximal smooth atlas with charts
$\varphi:U\to\varphi(U)\subseteq\mathbb R^n_+$,
where $\mathbb R^n_+=\{(x_1,\dotsc,x_n)\mid x_i\geq 0\}$.
The collection of points of $P$ which correspond by some
(and hence every) chart to points in $\mathbb R^n$ with
exactly $k$ coordinates equal to zero is a well defined subset
of $P$ called the \emph{$k$--corner of $P$} and it will be
denoted by
$P_k$. It has a structure of a smooth $(n-k)$--dimensional
manifold. The union $\partial P=P_1\cup P_2\cup\cdots\cup P_n$
is a closed subset which is a topological manifold and
$(P,\partial P)$ is a topological manifold with boundary
$\partial P$.

The following theorem can easily be derived from \cite[Theorem~1]{BH01}
by lifting everything to the universal covering, see
Proposition~\ref{P:cano}.

\begin{theorem}\label{T:6}
Let $M$ be a closed manifold, and suppose $X$ is a smooth vector field
which satisfies \MS\ and \LL.
\begin{enumerate}
\item\label{T:6i}
For any two rest points $\tilde x,\tilde y\in\tilde{\mathcal X}$ the set 
$\hat{\mathcal T}(\tilde x,\tilde y)$ admits a natural structure of a compact 
smooth manifold with corners, whose $k$--corner coincides 
with $\hat{\mathcal T}(\tilde x,\tilde y)_k$ from \eqref{Bk_defi}.
\item\label{T:6ii}
For every rest point $\tilde x\in\tilde{\mathcal X}$,
the set $\hat W_{\tilde x}^-$ admits a natural structure of a smooth
manifold with corners, whose $k$--corner coincides with
$(\hat W^-_{\tilde x})_k$ from~\eqref{Wk_defi}.
\item\label{T:6iii}
The function
$\hat i^-_{\tilde x}:\hat W^-_{\tilde x}\to\tilde M$
is smooth and proper, for all $\tilde x\in\tilde{\mathcal X}$.
\item\label{T:6iv} 
If $\omega$ is Lyapunov for $X$ and $h:\tilde M\to\mathbb R$ is a 
smooth function with $dh=\pi^*\omega$ then
the function $h\circ\hat i^-_{\tilde x}$ is
smooth and proper, for all $\tilde x\in\tilde{\mathcal X}$.
\end{enumerate}
\end{theorem}

As a first folklore application of Theorem~\ref{T:6} we will give a

\begin{proof}[Proof of Proposition~\ref{P:d2=0}]
Let $x,z\in\mathcal X$.
Theorem~\ref{T:6}\itemref{T:6i} implies 
$$
\sum_{y\in\mathcal X}
\sum_{\gamma_1\in\mathcal P_{x,y}}
\mathbb I_{x,y}(\gamma_1)
\cdot\mathbb I_{y,z}(\gamma_1^{-1}\gamma)
=0
$$
for all $\gamma\in\mathcal P_{x,z}$.
If $\eta\in \QX$ we can reorder sums and find
\begin{multline*}
\sum_{y\in\mathcal X}
\sum_{\gamma_1\in\mathcal P_{x,y}}
\mathbb I_{x,y}(\gamma_1)e^{\eta(\gamma_1)}
\sum_{\gamma_2\in\mathcal P_{y,z}}
\mathbb I_{y,z}(\gamma_2)e^{\eta(\gamma_2)}
\\=
\sum_{\gamma\in\mathcal P_{x,z}}
\Bigl(
\sum_{y\in\mathcal X}
\sum_{\gamma_1\in\mathcal P_{x,y}}
\mathbb I_{x,y}(\gamma_1)
\cdot\mathbb I_{y,z}(\gamma_1^{-1}\gamma)
\Bigr)
e^{\eta(\gamma)}
=0.
\end{multline*}
This implies $\delta_\eta^2=0$.
\end{proof}

As a second application of Theorem~\ref{T:6} we will give a

\begin{proof}[Proof of Proposition~\ref{P:inthom}]
We follow the approach in \cite{BH01}.
Let $\chi:\R\to[0,1]$ be smooth, and such that $\chi(t)=0$ for $t\leq0$ and
$\chi(t)=1$ for $t\geq1$. Choose a Lyapunov form $\omega$ for $X$.
For $y\in\mathcal X$ and $s\in\mathbb R$ define 
$\hat\chi^s_y:=\chi\circ(\hat h^\omega_y+s):\hat W_y^-\to[0,1]$.
Note that $\supp(\hat\chi_y^s)$ is compact in view of
Theorem~\ref{T:6}\itemref{T:6iv}.
Suppose $x\in\mathcal X$, $\alpha\in\Omega^*(M;\mathbb C)$, and $\eta\in \RX$.
Absolute convergence implies
$$
\Int_\eta(d_\eta\alpha)(x)=
\lim_{s\to\infty}\int_{\hat W_x^-}\hat\chi_x^s\cdot
e^{\hat h^\eta_x}\cdot(\hat i_x^-)^*d_\eta\alpha.
$$ 
Moreover,
$$
\hat\chi_x^s\cdot e^{\hat h^\eta_x}\cdot(\hat i_x^-)^*d_\eta\alpha
=d\bigl(\hat\chi_x^s\cdot e^{\hat h^\eta_x}\cdot(\hat i_x^-)^*\alpha\bigr)
-\bigl(\chi'\circ(\hat h_x^\omega+s)\bigr)\cdot e^{\hat h^\eta_x}
\cdot(\hat i_x^-)^*\omega\wedge\alpha.
$$
Since $\eta\in \RX$ and since $\chi'$ is bounded we have
$$
\lim_{s\to\infty}\int_{\hat W_x^-}\bigl(\chi'\circ(\hat h^\omega_x+s)\bigr)
\cdot e^{\hat h^\eta_x}\cdot (\hat i_x^-)^*\omega\wedge\alpha
=0.
$$
Using Theorem~\ref{T:6}\itemref{T:6ii} and 
Stokes' theorem for the compactly supported smooth form
$\hat\chi_x^s\cdot e^{\hat h_x^\eta}\cdot(\hat i_x^-)^*\alpha\in\Omega^*(\hat W_x^-;\mathbb C)$
we get
$$
\int_{\hat W_x^-}     
d\bigl(\hat\chi_x^s\cdot e^{\hat h^\eta_x}\cdot(\hat i_x^-)^*\alpha\bigr)
=\sum_{y\in\mathcal X}\sum_{\gamma\in\mathcal P_{x,y}}
\mathbb I_{x,y}(\gamma)e^{\eta(\gamma)}
\int_{\hat W_y^-}\hat\chi_y^{s+\omega(\gamma)}\cdot e^{\hat h^\eta_y}\cdot(\hat i_y^-)^*\alpha.
$$
Since $\eta\in \QX\cap \RX$ the form
$\mathbb I_{x,y}(\gamma)e^{\eta(\gamma)}\cdot e^{\hat h^\eta_y}\cdot(\hat i_y^-)^*\alpha$
is absolutely integrable on $\mathcal P_{x,y}\times\hat W_y^-$. Hence we may
interchange limits and find
\begin{multline*}
\lim_{s\to\infty}\int_{\hat W_x^-}
d\bigl(\hat\chi_x^s\cdot e^{\hat h^\eta_x}\cdot(\hat i_x^-)^*\alpha\bigr)
\\
=\sum_{y\in\mathcal X}\sum_{\gamma\in\mathcal P_{x,y}}
\mathbb I_{x,y}(\gamma)e^{\eta(\gamma)}
\lim_{s\to\infty}\int_{\hat W_y^-}
\hat\chi_x^{s+\omega(\gamma)}\cdot e^{\hat h^\eta_y}\cdot(\hat i_y^-)^*\alpha
\\               
=\sum_{y\in\mathcal X}L(\mathbb I_{x,y})(\eta)\cdot
\Int_\eta(\alpha)(y)
=\delta_\eta(\Int_\eta(\alpha))(x).
\end{multline*}  
We conclude      
$\Int_\eta(d_\eta\alpha)(x)=\delta_\eta(\Int_\eta(\alpha))(x)$.
\end{proof}

We close this section with a lemma which immediately implies Proposition~\ref{P:onto}.

\begin{lemma}\label{L:inte}
Suppose $\eta\in \RX$, $x\in\mathcal X$, and let $\epsilon>0$.
Then there exists $\alpha\in\Omega^*(M;\mathbb C)$ so that
$|\Int_\eta(\alpha)(y)-\delta_{x,y}|\leq\epsilon$, for all $y\in\mathcal X$.
\end{lemma}

\begin{proof}
We follow the approach in \cite{BH01}. Let $U$ be a neighborhood of $x$
on which $X$ has canonical form \eqref{E:4}.
Let $B\subseteq W_x^-$ denote the connected component of $W_x^-\cap U$
containing $x$.
Choose $\alpha\in\Omega^*(M;\mathbb C)$ with $\supp(\alpha)\subseteq U$ and such
that
$\int_Be^{h^\eta_x}(i_x^-)^*\alpha=1$.
For every $y\in\mathcal X$ choose a compact $K_y\subseteq W_y^-$ such that
$\int_{W_y^-\setminus K_y}|e^{h^\eta_y}(i_y^-)^*\alpha|\leq\epsilon$.
Assume $B\subseteq K_x$.
By multiplying $\alpha$ with a bump function which is $1$ on
$B$ and whose support is sufficiently concentrated around $B$,
we may in addition assume $\supp(\alpha)\cap(K_x\setminus B)=\emptyset$, and
$\supp(\alpha)\cap K_y=\emptyset$ for all $x\neq y\in\mathcal X$.
Then
\[
\bigl|\Int_\eta(\alpha)(y)-\delta_{x,y}\bigr|
=\Bigl|\int_{W_y^-\setminus K_y}e^{h^\eta_y}(i_y^-)^*\alpha\Bigr|
\leq\int_{W_y^-\setminus K_y}|e^{h^\eta_y}(i_y^-)^*\alpha|\leq\epsilon.
\qedhere
\]
\end{proof}

\subsection{Proof of the first part of Theorem~\ref{T:int}}\label{SS:Tinti}

Suppose $X$ satisfies \MS\ and \LL.
Let $\Gamma:=\pi_1(M)$ denote the fundamental group acting from the left on
the universal covering $\pi:\tilde M\to M$ in the usual manner.
Equip $\mathbb C^{\mathcal X}$ with a norm.
Equip $\mathbb A:=\eend(\mathbb C^{\mathcal X})$
with the corresponding operator norm.
Choose a Lyapunov form $\omega$ for $X$.
Let $N$ denote the vector space of maps $a:\Gamma\to\mathbb A$
for which 
$\{\gamma\in\Gamma\mid -\omega(\gamma)\leq K, a(\gamma)\neq0\}$ 
is finite, for all $K\in\mathbb R$.
Equipped with the convolution product $N$ becomes an algebra with unit.
Let $L^1:=L^1(\Gamma;\mathbb A)$ denote the Banach space
of functions $a:\Gamma\to\mathbb A$ for which
$\|a\|_{L^1}:=\sum_{\gamma\in\Gamma}\|a(\gamma)\|<\infty$.
Recall that the convolution product makes $L^1$ a Banach algebra
with unit.

\begin{lemma}\label{L:NL1}
Let $I,a\in N$. Assume $\|1-a\|_{L^1}<1$, $I*a\in L^1$, and assume that 
$a(\gamma)\neq0$ implies $-\omega(\gamma)\geq0$. Then $I\in L^1$.
\end{lemma}

\begin{proof}
Since $L^1$ is a Banach algebra
$\|1-a\|_{L^1}<1$ implies that $a\in L^1$ is invertible.
Clearly it suffices to show $(I*a)*a^{-1}=I*(a*a^{-1})$.
That is, for fixed $\rho\in\Gamma$, we have to show
\begin{equation}\label{E:xyz}
\sum_{\sigma\in\Gamma}\sum_{\tau\in\Gamma}
I(\sigma)a(\sigma^{-1}\tau)a^{-1}(\tau^{-1}\rho)
=
\sum_{\tau\in\Gamma}\sum_{\sigma\in\Gamma}
I(\sigma)a(\sigma^{-1}\tau)a^{-1}(\tau^{-1}\rho).
\end{equation}
Using $a^{-1}=\sum_{k=0}^\infty(1-a)^k$ it is not difficult to show that 
$a^{-1}(\gamma)\neq0$ implies $-\omega(\gamma)\geq0$. Using $I,a\in N$ we thus conclude that
$$
\bigl\{\sigma\in\Gamma\bigm|
\exists\tau\in\Gamma:I(\sigma)a(\sigma^{-1}\tau)a^{-1}(\tau^{-1}\rho)\neq0\bigr\}
$$
is finite. Equation~\eqref{E:xyz} follows immediately.
\end{proof}

Let us now turn to the proof of $\RX\subseteq\QX$. 
Let $\eta\in \RX$. Choose a lift $s(x)\in\tilde{\mathcal X}$ 
for every zero $x\in\mathcal X$, i.e.\ $\pi(s(x))=x$. 
For $x,y\in\mathcal X$ and $\gamma\in\Gamma$ let $\rho_{x,y}^\gamma\in\mathcal P_{x,y}$ denote
the homotopy class of paths determined by the lifts $s(x)$ and $\gamma\cdot s(y)$.
Moreover, set
\begin{equation}\label{E:15}
I_{x,y}(\gamma):=\mathbb I_{x,y}(\rho_{x,y}^\gamma)
\cdot e^{\eta(\rho_{x,y}^\gamma)}.
\end{equation}
We write $I:\Gamma\to\mathbb A$ for the matrix valued map defined by
\eqref{E:15}. Note that $\eta\in \QX$ iff $I\in L^1$.
It thus suffices to construct $a:\Gamma\to\mathbb A$ for which
Lemma~\ref{L:NL1} is applicable.
Note that $I\in N$ in view of Proposition~\ref{P:Nov}.

In order to construct $a$ choose a smooth function
$\chi:\tilde M\to[0,1]$ so that $\chi=1$ in a
neighborhood of $s(\mathcal X)$, so that $\supp(\chi)$ is compact, and so that
$\supp(\chi)\cap\supp(\gamma^*\chi)=\emptyset$ for all non-trivial
$\gamma\in\Gamma$.
For $x\in\mathcal X$ and $\gamma\in\Gamma$ define a function
$\hat\chi^\gamma_x:\hat W_x^-\to[0,1]$ by 
$\hat\chi_x^\gamma:=(\gamma^{-1})^*\chi\circ\hat i_{s(x)}^-$.
Note that $\supp(\hat\chi_x^\gamma)$ is compact in view of
Theorem~\ref{T:6}\itemref{T:6iii}.
Possibly shrinking the support of $\chi$ we may assume that
$\hat\chi_x^\gamma\neq0$ implies $-\omega(\gamma)\geq0$.

The construction of $a$ will also depend on the choice of 
$\beta_x\in\Omega^*(M;\mathbb C)$, $x\in\mathcal X$, 
which will be specified below. For $x,y\in\mathcal X$ and $\gamma\in\Gamma$
define
\begin{equation}\label{E:aee}
a_{x,y}(\gamma)
:=\int_{\hat W_x^-}\hat\chi_x^\gamma\cdot e^{\hat h^\eta_x}\cdot(\hat i_x^-)^*\beta_y.
\end{equation}
We write $a:\Gamma\to\mathbb A$ for the matrix valued function defined by
\eqref{E:aee}. Note that $a\in N$.
Moreover, $a(\gamma)\neq0$ implies $-\omega(\gamma)\geq0$.

We will choose $\beta_x$ so that its support is concentrated near $x$.
More precisely, we assume $\supp(d\chi)\cap\supp(\pi^*\beta_x)=\emptyset$ 
for all $x\in\mathcal X$.
Clearly we may also assume that $a_{x,y}(e)=\delta_{x,y}$, i.e.\
$a(e)=1$. Note that the mutual disjointness
of $\supp(\chi_x^\gamma)$, $\gamma\in\Gamma$, implies
$$
\sum_{e\neq\gamma\in\Gamma}|a_{x,y}(\gamma)|
\leq
\sum_{e\neq\gamma\in\Gamma}
\int_{\supp(\hat\chi_x^\gamma)}|e^{\hat h^\eta_x}\cdot(\hat i_x^-)^*\beta_y|
\leq
\int_{\hat W_x^-\setminus\supp(\hat\chi_x^e)}|e^{\hat h^\eta_x}\cdot(\hat i_x^-)^*\beta_y|.
$$
Using $\eta\in \RX$ and arguing as in the proof of Lemma~\ref{L:inte}, 
we may therefore assume that, given $\epsilon>0$, the $\beta_x$ are chosen so that
$\sum_{e\neq\gamma\in\Gamma}|a_{x,y}(\gamma)|<\epsilon$ for all
$x,y\in\mathcal X$.
Obviously this implies $\|1-a\|_{L^1}<1$.

Using $\supp(d\chi)\cap\supp(\pi^*\beta_y)=\emptyset$ and applying
Stokes' theorem for the compactly supported form
$\hat\chi_x^\gamma\cdot e^{\hat h^\eta_x}\cdot(\hat i_x^-)^*\beta_y
\in\Omega^*(\hat W_x^-;\mathbb C)$, see Theorem~\ref{T:6}\itemref{T:6ii}, we find
\begin{multline*}
\int_{\hat W_x^-}\hat\chi_x^\gamma\cdot e^{\hat h_x^\eta}\cdot(\hat i_x^-)^*d_\eta\beta_y
=\int_{\hat W_x^-}d\bigl(\hat\chi_x^\gamma\cdot e^{\hat h_x^\eta}\cdot(\hat i_x^-)^*\beta_y\bigr)
\\
=\sum_{z\in\mathcal X}\sum_{\sigma\in\Gamma}
\mathbb I_{x,z}(\alpha_{xz}^\sigma)e^{\eta(\alpha_{xz}^\sigma)}
\int_{\hat W_z^-}\hat\chi^{\sigma^{-1}\gamma}_z\cdot e^{\hat h_z^\eta}
\cdot(\hat i_z^-)^*\beta_y
=(I*a)_{x,y}(\gamma).
\end{multline*}
Since $\eta\in \RX$ we therefore get
$$
\sum_{\gamma\in\Gamma}|(I*a)_{x,y}(\gamma)|
\leq\sum_{\gamma\in\Gamma}\int_{\supp(\hat\chi_x^\gamma)}
|e^{\hat h^\eta_x}\cdot(\hat i_x^-)^*d_\eta\beta_y|
\leq\int_{\hat W_x^-}|e^{\hat h_x^\eta}\cdot(\hat i_x^-)d_\eta\beta_y|<\infty
$$
for all $x,y\in\mathcal X$. We conclude $\|I*a\|_{L^1}<\infty$.
Hence we can apply Lemma~\ref{L:NL1}, obtain $I\in L^1$ and thus $\eta\in \QX$.
This completes the proof of $\RX\subseteq \QX$.

\subsection{Proof of the second part of Theorem~\ref{T:int}}\label{SS:Tintii}

We will start with a lemma whose first part, when applied to the
eigen spaces of $B_\eta$, implies Proposition~\ref{P:hodge}.

\begin{lemma}\label{L:finHodge}
Let $C^*$ be a finite dimensional graded complex over $\mathbb C$ with differential $d$.
Let $b$ be a non-degenerate graded bilinear form on $C^*$.
Let $d^t$ denote the formal transpose of $d$, i.e.\
$b(dv,w)=b(v,d^tw)$ for all $v,w\in C^*$.
Set $B:=dd^t+d^td$ and suppose $\ker B=0$. Then $C^*=\img d\oplus\img d^t$,
and this decomposition is orthogonal with respect to $b$. Particularly, the cohomology
of $C^*$ vanishes. For its torsion, with respect to the equivalence class of graded bases 
determined by $b$, we have
$$
\tau(C^*,b)^2=\prod_q(\det B^q)^{(-1)^{q+1}q}
$$
where $B^q:C^q\to C^q$ denotes the part of $B$ acting in degree $q$.
\end{lemma}

\begin{proof}
Clearly $\img d\subseteq(\ker d^t)^\perp$, and hence
$\img d=(\ker d^t)^\perp$ since $C^*$ is finite dimensional. 
Similarly we get $\img d^t=(\ker d)^\perp$.
Therefore
$$
(\img d+\img d^t)^\perp
=(\img d)^\perp\cap(\img d^t)^\perp
=\ker d^t\cap\ker d
\subseteq\ker B=0,
$$
and thus $C^*=\img d+\img d^t$. Moreover, since 
$\img d^t\subseteq(\ker d)^\perp\subseteq(\img d)^\perp$
this decomposition is orthogonal. The cohomology vanishes
for we have $\ker d=(\img d^t)^\perp=\img d$.
Using	
\begin{eqnarray*}
\det B^q
&=&
\det(d^td|_{\img d^t\cap C^q})\cdot\det(dd^t|_{\img d\cap C^q})
\\&=&
\det(d^td|_{\img d^t\cap C^q})\cdot\det(d^td|_{\img d^t\cap C^{q-1}})
\end{eqnarray*}
a trivial telescoping calculation shows
\[
\prod_q(\det B^q)^{(-1)^{q+1}q}
=\prod_q\bigl(\det d^td|_{\img d^t\cap C^q}\bigr)^{(-1)^q}
=\tau(C^*,b)^2.
\qedhere
\]
\end{proof}

Let us next prove that $\mathring\RX\cap\SX$ is an analytic subset
of $\mathring\RX$. Let $\eta_0\in\mathring\RX$. Choose a simple 
closed curve $K$ around $0\in\mathbb C$
which avoids the spectrum of $B_{\eta_0}$. Let $U$
be an open neighborhood of $\eta_0$ so that $K$ avoids the spectrum of every 
$B_\eta$, $\eta\in U$. Assume $U$ is connected and $U\subseteq\mathring \RX$.
Let $E^*_\eta(K)$ denote the image of the spectral
projection associated with $K$, i.e.\ $E^*_\eta(K)$ is the sum of
all eigen spaces of $B_\eta$ corresponding to eigen values contained in the interior of
$K$. Since the spectral projection depends holomorphically 
on $\eta$, wee see that $E^*_\eta(K)$ is a holomorphic family of 
finite dimensional complexes parametrized by $\eta\in U$. From Proposition~\ref{P:hodge} we see
that the inclusion $E^*_\eta(K)\to\Omega^*_\eta(M;\mathbb C)$ is a quasi
isomorphism for all $\eta\in U$.

Consider the restriction of the integration
$\Int_\eta:E^*_\eta(K)\to C^*_\eta(X;\mathbb C)$,
and let $\mathcal C_\eta^*(K)$ denote its mapping cone. 
More precisely, as graded vector space $\mathcal C_\eta^*(K)=
C^{*-1}_\eta(X;\mathbb C)\oplus E^*_\eta(K)$, and the differential is given by
$(f,\alpha)\mapsto(-\delta_\eta f+\Int_\eta\alpha,d_\eta\alpha)$.
This is a family
of finite dimensional complexes, holomorphically parametrized by $\eta\in U$. Note that the 
dimension of $\mathcal C^*_\eta(K)$ is even, $\dim\mathcal C^*_\eta(K)=2k$. 
Possibly shrinking $U$ we may assume that we have a base
$\{v^1_\eta,\dotsc,v^{2k}_\eta\}$ of $\mathcal C^*_\eta(K)$ holomorphically parametrized by
$\eta\in U$. Let $f^1_\eta,\dotsc,f^N_\eta\in\mathbb C$ denote the 
$k\times k$--minors of the differential of $\mathcal C_\eta^*(K)$ 
with respect to this base, 
$N=\left(\begin{smallmatrix}2k\\ k\end{smallmatrix}\right)^2$. This provides $N$ holomorphic
functions $f^i:U\to\mathbb C$, $1\leq i\leq N$. 
For $\eta\in U$ the integration will induce
an isomorphism in cohomology iff $\mathcal C^*_\eta(K)$ is acyclic.
Moreover, $\mathcal C^*_\eta(K)$ is acyclic iff its differential has rank $k$.
This in turn is equivalent to $f^i(\eta)\neq0$ for some $1\leq i\leq N$.
We conclude $\SX\cap U=\{\eta\in U\mid f^i(\eta)=0, 1\leq i\leq N\}$.
Hence $\mathring\RX\cap\SX$ is an analytic subset of $\mathring\RX$.

Suppose $\omega$ is a Lyapunov form for $X$, and let $\eta\in\RX$.
Recall that we have an integration homomorphism
\begin{equation}\label{E:368}
\Int_{\eta+t\omega}:\Omega^*_{\eta+t\omega}(M;\mathbb C)\to
C^*_{\eta+t\omega}(X;\mathbb C)
\end{equation}
for $t\geq0$, and that $\eta+t\omega\in\mathring\RX$ for $t>0$, see
Proposition~\ref{P:basic:R}. We have to show that \eqref{E:368}
induces an isomorphism in cohomology for sufficiently large $t$.
In view of the gauge invariance, see \eqref{CD:gauge},
we may assume that $\eta$ vanishes in a neighborhood of $\mathcal X$, and
that there exists a Riemannian metric $g$, such that $\omega=-g(X,\cdot)$
as in Proposition~\ref{P:cano}.

Consider the one parameter family of complexes
$\Omega^*_{\eta+t\omega}(M;\mathbb C)$. Let $\Delta_{\eta+t\omega}$
denote the corresponding Laplacians with respect to the standard Hermitian
structure on $\Omega^*(M;\mathbb C)$. Witten--Helffer--Sj\"ostrand theory
\cite{BH01} tells that as $t\to\infty$ the spectrum of
$\Delta_{\eta+t\omega}$ develops a gap, providing a canonic orthogonal
decomposition
$$
\Omega^*_{\eta+t\omega}(M;\mathbb C)
=\Omega^*_{\eta+t\omega,\sm}(M;\mathbb C)
\oplus\Omega^*_{\eta+t\omega,\la}(M;\mathbb C).
$$
Moreover, for sufficiently large $t$ the restriction of the integration
$$
\Int_{\eta+t\omega}:\Omega^*_{\eta+t\omega,\sm}(M;\mathbb C)\to
C^*_{\eta+t\omega}(X;\mathbb C)
$$
is an isomorphism \cite{BH01}. It follows that
\eqref{E:368} induces an isomorphism in cohomology, and hence
$\eta+t\omega\in\mathring\RX\setminus\Sigma$ for
sufficiently large $t$.

\section{Proof of Theorem~\ref{T:tor}}\label{S:Ttor}

The proof of Theorem~\ref{T:tor} is based a result of Bismut--Zhang
\cite{BZ92} and formula of Hutchings--Lee \cite{HL99} and Pajitnov \cite{P03}.
The Bismut--Zhang theorem implies that Theorem~\ref{T:tor} holds
for Morse--Smale vector fields, see section~\ref{SS:BZ}.
The Hutchings--Lee formula permits to establish an anomaly formula
in $X$ for $(T\Int_\eta^X)^2$, see Proposition~\ref{P:inp_anom} in section~\ref{SS:HLP}.
Putting this together we will obtain Theorem~\ref{T:tor}, see
section~\ref{SS:Pthmint}.

\subsection{Proof of Proposition~\ref{P:Tint}}\label{SS:PTint}

Let us first show $(T\Int_{\bar\eta})^2=\overline{(T\Int_\eta)^2}$.
Clearly we have $R(\bar\eta,X,g)=\overline{R(\eta,X,g)}$.
Note that complex conjugation on $\Omega^*(M;\mathbb C)$
intertwines $d_\eta$ with $d_{\bar\eta}$, $d^t_\eta$
with $d^t_{\bar\eta}$, and $B_\eta$ with $B_{\bar\eta}$.
Therefor the spectrum of $B_\eta$ is conjugate to the spectrum
of $B_{\bar\eta}$. It follows that 
$(T^\an_{\bar\eta,g})^2=\overline{(T^\an_{\eta,g})^2}$.
Moreover, complex conjugation restricts to an anti-linear isomorphism of complexes
$E^*_\eta(0)\simeq E^*_{\bar\eta}(0)$ which is easily seen to intertwine
the equivalence class of bases determined by $b$.
Complex conjugation also defines an anti-linear isomorphism of complexes
$C^*_\eta(X;\mathbb C)\simeq C^*_{\bar\eta}(X;\mathbb C)$ which intertwines the
equivalence class of bases determined by the indicator functions.
These isomorphisms intertwine $\Int_\eta$ with $\Int_{\bar\eta}$.
Hence they provide an anti-linear isomorphism of mapping cones, and therefore
$\pm T(\Int_{\bar\eta}|_{E^*_{\bar\eta}(0)})=\overline{\pm T(\Int_\eta|_{E^*_\eta(0)})}$.
Putting everything together we find
$(T\Int_{\bar\eta})^2=\overline{(T\Int_\eta)^2}$.

Let us next show that $(T\Int_\eta)^2$ depends holomorphically on
$\eta\in\mathring \RX\setminus\SX$. Let $\eta_0\in\mathring
\RX\setminus\SX$. As in the proof of Theorem~\ref{T:int} in
section~\ref{SS:Tintii} let $U$ be a connected open neighborhood of $\eta_0$
so that $K$ avoids the spectrum of $B_\eta$ for all $\eta\in U$.
Assume $U\subseteq\mathring \RX\setminus\SX$.
For $\eta\in U$ let us write 
$\prod_q(\det{}^K B^q_\eta)^{(-1)^{q+1}q}$
for the zeta regularized product of eigen values of $B_\eta$ not contained in the
interior of $K$. This depends holomorphically on $\eta\in U$.
Let us write $\mathcal C^*_\eta(K)$ for the mapping cone of 
$\Int_\eta:E^*_\eta(K)\to C^*_\eta(X;\mathbb C)$. This is a finite
dimensional family of complexes holomorphically parametrized by
$\eta\in U$, see section~\ref{SS:Tintii}.
Note that these complexes are acyclic since $U\cap\SX=\emptyset$.
We equip $\mathcal C^*_\eta(K)$ with the basis determined by the restriction
of the bilinear form $b$ and the indicator functions in $C^*_\eta(X;\mathbb C)$. 
These equivalence classes of bases depend holomorphically on $\eta\in U$.
Hence the torsion 
$(T(\Int_\eta|_{E^*_\eta(K)}))^2=(T\mathcal C^*_\eta(K))^2$ depends holomorphically on
$\eta\in U$. Using Lemma~\ref{L:finHodge} it is easy to see that
$$
(T(\Int_\eta|_{E_\eta^*(0)}))^2\cdot\prod_q(\det{}'B_\eta^q)^{(-1)^{q+1}q}
=(T(\Int_\eta|_{E_\eta^*(K)}))^2\cdot\prod_q(\det{}^KB_\eta^q)^{(-1)^{q+1}q}.
$$
Hence $(T\Int_\eta)^2$ depends holomorphically on $\eta$ too.

Similarly, using \eqref{E:cont:I} and \eqref{E:cont:R} 
one shows that $\lim_{t\to0^+}(T\Int_{\eta+t\omega})^2=(T\Int_\eta)^2$
for a Lyapunov form $\omega$ and $\eta\in\RX\setminus\Sigma$.

Next we will show that $(T\Int_\eta)^2$ does not depend on $g$. 
For real valued $\eta\in\mathcal Z^1(M;\mathbb R)\cap(\RX\setminus\Sigma)$
the operator $B_\eta$ coincides with the Laplacian associated with $g$ and
$\eta$, and hence $T^\an_{\eta,g}$ coincides with the Ray--Singer torsion.
For two Riemannian metrics $g_1$ and $g_2$ on $M$,
the anomaly formula in \cite[Theorem~0.1]{BZ92} then implies
$$
\log\frac{(T(\Int_\eta|_{E_{\eta,g_2}^*(0)}))^2\cdot(T^\an_{\eta,g_2})^2}
{(T(\Int_\eta|_{E_{\eta,g_1}^*(0)}))^2\cdot(T^\an_{\eta,g_1})^2}
=2\int_M\eta\wedge\cs(g_1,g_2).
$$
Together with \eqref{E:rg1g2}
this yields $(T\Int_{\eta,g_1})^2=(T\Int_{\eta,g_2})^2$.
Since both sides depend holomorphically on $\eta$, see
Proposition~\ref{P:Tint}, this relation is true for
$\eta\in\mathring\RX\setminus\Sigma$ too. In view of \eqref{E:cont:int}
it continues to hold for $\eta\in\RX\setminus\Sigma$.

Let us finally turn to the gauge invariance. 
Again, for real $\eta\in\mathcal Z^1(M;\mathbb R)\cap(\RX\setminus\Sigma)$ and 
real $h\in C^\infty(M;\mathbb R)$ the anomaly formula in
\cite[Theorem~0.1.]{BZ92} implies
$$
\log\frac{(T(\Int_\eta|_{E_{\eta+dh,g}^*(0)}))^2\cdot(T^\an_{\eta+dh,g})^2}
{(T(\Int_\eta|_{E_{\eta,g}^*(0)}))^2\cdot(T^\an_{\eta,g})^2}
=2\biggl(-\int_Mh\ec(g)+\sum_{x\in\mathcal X}\IND(x)h(x)\biggr).
$$
Together with \eqref{E:Rdh} this implies
$(T\Int_{\eta+dh})^2=(T\Int_\eta)^2$.
Since both sides depend holomorphically on $\eta$ and $h$, see
Proposition~\ref{P:Tint},
this relation continues to hold for $\eta\in\mathring\RX\setminus\Sigma$
and $h\in C^\infty(M;\mathbb C)$.
In view of \eqref{E:cont:int} it remains true for
$\eta\in\RX\setminus\Sigma$.
This completes the proof of Proposition~\ref{P:Tint}.

\subsection{The Bismut--Zhang theorem}\label{SS:BZ}

Suppose our vector field is of the form $X=-\grad_{g_0}f$ for some
Riemannian metric $g_0$ on $M$ and a Morse function $f:M\to\mathbb R$.
Then $df$ is Lyapunov for $X$, hence $X$ satisfies \LL.
There are no closed trajectories. Hence $X$ satisfies \NCT,
$\PX=\mathcal Z^1(M;\mathbb C)$ and $e^{L(h_*\mathbb P)(\eta)}=1$.
In view of Theorem~\ref{T:6}\itemref{T:6iv}
the completion of the unstable manifolds are compact, hence
$\RX=\QX=\mathcal Z^1(M;\mathbb C)$.
It is well known that $\Sigma=\emptyset$, i.e.\ 
the integration induces an isomorphism for all $\eta$.
A theorem of Bismut--Zhang \cite[Theorem~0.2]{BZ92} tells that 
in this case
\begin{equation}\label{E:BZ}
(T\Int_\eta)^2=1
\end{equation}
for all $\eta\in\mathcal Z^1(M;\mathbb R)$.
Since $(T\Int_\eta)^2$ depends holomorphically on $\eta$, see
Proposition~\ref{P:Tint}, the relation \eqref{E:BZ} continues to hold for all 
$\eta\in\mathcal Z^1(M;\mathbb C)$. To make a long story short, Theorem~\ref{T:tor}
is true for vector fields of the form $X=-\grad_{g_0}f$.

\subsection{An anomaly formula}\label{SS:anom}

Consider the bordism $W:=M\times[-1,1]$.
Set $\partial_\pm W:=M\times\{\pm1\}$.
Let $Y$ be a vector field on $W$.
Assume that there are vector fields $X_\pm$ on $M$ so that
$Y(z,s)=X_+(z)+(s-1)\partial/\partial s$
in a neighborhood of $\partial_+W$ and so that
$Y(z,s)=X_-(z)+(-s-1)\partial/\partial s$ in a neighborhood of $\partial_-W$.
Particularly, $Y$ is tangential to $\partial W$.
Moreover, assume that $ds(Y)<0$ on $M\times(-1,1)$.
Particularly, there are no zeros or closed trajectories of $Y$
contained in the interior of $W$.
Let $\mathcal X_\pm$ denote the zeros of $X_\pm$.
For $x\in\mathcal X_-$ we have $\ind_Y(x)=\ind_{X_-}(x)$, but note that
for $x\in\mathcal X_+$ we have $\ind_Y(x)=\ind_{X_+}(x)+1$.
We choose the orientations of the unstable manifolds of $Y$ so that 
$W_{Y,x}^-=W_{X_-,x}^-$ is orientation preserving for $x\in\mathcal X_-$,
and so that $\partial W^-_{Y,x}=W^-_{X_+,x}$ is orientation reversing
for $x\in\mathcal X_+$.

Suppose $Y$ satisfies \MS\ and \LL. Note that this implies that $X_\pm$ satisfy
\MS\ and \LL\ too. Then Proposition~\ref{P:Nov}, Theorem~\ref{T:6}\itemref{T:6i}
and Proposition~\ref{P:d2=0} continue to hold for $Y$. Hence we get a complex
$C^*_{\tilde\eta}(Y;\mathbb C)$ for all ${\tilde\eta}\in\QX^Y$. Note that for 
$\tilde\eta\in\QX^Y$ we have $\eta_\pm:=\iota_\pm^*\tilde\eta\in\QX^{X_\pm}$,
where $\iota_\pm:M\to\partial_\pm W$, $\iota_\pm(z)=(z,\pm1)$.
Clearly
\begin{equation}\label{E:scy}
C^*_{\tilde\eta}(Y;\mathbb C)
=C^{*-1}_{\eta_+}(X_+;\mathbb C)\oplus C^*_{\eta_-}(X_-;\mathbb C),
\qquad
\delta^Y_{\tilde\eta}=
\left(
\begin{smallmatrix}
-\delta^{X_+}_{\eta_+} & u^Y_\eta
\\
0 & \delta^{X_-}_{\eta_-}
\end{smallmatrix}
\right)
\end{equation}
for some 
\begin{equation}\label{E:432}
u_{\tilde\eta}^Y:C^*_{\eta_-}(X_-;\mathbb C)\to C^*_{\eta_+}(X_+;\mathbb C).
\end{equation}
>From $(\delta^Y_{\tilde\eta})^2=0$ we see that \eqref{E:432}
is a homomorphism of complexes.

Theorem~\ref{T:6}\itemref{T:6ii} needs a minor adjustment in the case with
boundary. More precisely, for $x\in\mathcal X_+$ the 
completion of the unstable manifold $W_{x}^-$ has additional boundary parts
stemming from the fact that $W_x^-$ intersects $\partial_+W$ transversally.
For $\tilde\eta\in\RX^Y$ we get a linear mapping
$\Int_{\tilde\eta}^Y:\Omega^*(W;\mathbb C)\to C^*_{\tilde\eta}(Y;\mathbb C)$
satisfying
\begin{equation}\label{E:654}
\Int^Y_{\tilde\eta}\circ d_{\tilde\eta}
=\delta_{\tilde\eta}^Y\circ\Int^Y_{\tilde\eta}
-(\iota_+)_*\circ\Int_{\eta_+}^{X_+}\circ\iota^*_+.
\end{equation}
Here $\iota^*_+:\Omega^*_{\tilde\eta}(W;\mathbb C)\to\Omega^*_{\eta_+}(M;\mathbb C)$ is the pull
back of forms, and $(\iota_+)_*:C^*_{\eta_+}(X_+;\mathbb C)\to
C^{*+1}_{\tilde\eta}(Y;\mathbb C)$ is the obvious inclusion stemming from
$\mathcal X_+\subseteq\mathcal Y$. But note that while $\iota_+^*$ is a
homomorphism of complexes, we have 
$(\iota_+)_*\circ\delta_{\eta_+}^{X_+}+\delta_{\tilde\eta}^Y\circ(\iota_+)_*=0$.
Moreover, note that $\tilde\eta\in\RX^Y$ implies $\eta_\pm\in\RX^{X_\pm}$.
For $\eta_-$ this is trivial. For $\eta_+$ it follows from 
$W^-_{Y,x}\supseteq W^-_{X_+,x}\times(1-\epsilon,1]$ for some $\epsilon>0$.
Moreover, $\RX^Y\subseteq\QX^Y$, cf.\ Theorem~\ref{T:int}. So \eqref{E:654} indeed makes sense 
for $\tilde\eta\in\RX^Y$.
Splitting $\Int_{\tilde\eta}^Y$ according to \eqref{E:scy}
we find $\Int_{\tilde\eta}^Y=(h^Y_{\tilde\eta},\Int_{\eta_-}^{X_-}\circ\iota_-^*)$
for some 
$$
h_{\tilde\eta}^Y:\Omega^*_{\tilde\eta}(W;\mathbb C)
\to C^{*-1}_{\eta_+}(X_+;\mathbb C),
$$
and \eqref{E:654} tells that for all $\tilde\eta\in\RX^Y$
\begin{equation}\label{E:321}
h^Y_{\tilde\eta}\circ d_{\tilde\eta}
=-\delta^{X_+}_{\eta_+}\circ h^Y_{\tilde\eta}
+u^Y_{\tilde\eta}\circ\Int_{\eta_-}^{X_-}\circ\iota^*_-
-\Int^{X_+}_{\eta_+}\circ\iota^*_+.
\end{equation}

Let $p:W\to M$ denote the projection.
For $\eta\in\mathcal Z^1(M;\mathbb C)$ we write
$u^Y_\eta:=u^Y_{p^*\eta}$ and $h^Y_\eta:=h^Y_{p^*\eta}\circ p^*$.
Then $u^Y_\eta:C^*_\eta(X_-;\mathbb C)\to C^*_\eta(X_+;\mathbb C)$
is a homomorphism of complexes, and
$h^Y_\eta$ is a homotopy between $u^Y_\eta\circ\Int_\eta^{X_-}$ 
and $\Int_\eta^{X_+}$.

\begin{proposition}\label{P:anom}
Let $Y$ be a vector field on $W=M\times[-1,1]$ as above. Suppose
$\eta\in(\RX^{X_-}\setminus\Sigma^{X_-})\cap(\RX^{X_+}\setminus\Sigma^{X_+})$
and assume $p^*\eta\in\RX^Y$. Then $u^Y_\eta:C^*_\eta(X_-;\mathbb C)\to
C^*_\eta(X_+;\mathbb C)$ is a quasi isomorphism, and
\begin{equation}\label{E:anom}
\frac{(T\Int_{\eta}^{X_+})^2}{(T\Int_{\eta}^{X_-})^2}
=(Tu^Y_\eta)^2\cdot(e^{-\eta(\cs(X_-,X_+))})^2.
\end{equation}
Here the torsion $\pm Tu^Y_\eta$ is computed with respect to the base determined
by the indicator functions on $\mathcal X_\pm$.
\end{proposition}

\begin{proof}
>From the discussion above we know that $\Int_\eta^{X_+}$ is homotopic
to $u_\eta^Y\circ\Int^{X_-}_\eta$. Hence $u^Y_\eta$ is
a quasi isomorphism and
$$
\frac{\pm T(\Int^{X_+}_\eta|_{E^*_{\eta}(0)})}
{\pm T(\Int^{X_-}_\eta|_{E^*_{\eta}(0)})}
=\pm Tu_\eta^Y.
$$
Together with \eqref{EE:38} this yields \eqref{E:anom}.
\end{proof}

\subsection{Hutchings--Lee formula}\label{SS:HLP}

Let $X$ be a vector field which satisfies \MS\ and \LL.
Let $\Gamma:=\img\bigl(\pi_1(M)\to H_1(M;\mathbb R)\bigr)$.
Let $\omega\in\Omega^1(M;\mathbb R)$ be Lyapunov for $X$ and such that
$\omega:\Gamma\to\mathbb R$ is injective. 
Note that such Lyapunov forms exist in view of Proposition~\ref{P:cone}.
Let $\Lambda_\omega$
denote the corresponding Novikov field consisting of all functions
$\lambda:\Gamma\to\mathbb C$ for which 
$\{\gamma\in\Gamma\mid \lambda(\gamma)\neq0, -\omega(\gamma)\leq K\}$
is finite for all $K\in\mathbb R$, equipped with the convolution product.
Let us write $\Lambda_\omega^+$ for the subring of functions
$\lambda$ for which $\lambda(\gamma)\neq0$ implies $-\omega(\gamma)>0$.

The vector field $X$ gives rise to a Novikov complex
$C^*(X;\Lambda_\omega)$. This complex can be described as follows.
Let $\pi:\tilde M\to M$ denote the covering corresponding to the 
kernel of $\pi_1(M)\to H_1(M;\mathbb R)$. This is a principal
$\Gamma$--covering. Let $\tilde{\mathcal X}:=\pi^{-1}(\mathcal X)$
denote the zero set of the vector field $\tilde X:=\pi^*X$.
Choose a function $h:\tilde M\to\mathbb R$ such that $dh=\pi^*\omega$.
Now $C^*(X;\Lambda_\omega)$ is the space of all functions 
$c:\tilde{\mathcal X}\to\mathbb C$ for which
$\{\tilde x\in\tilde{\mathcal X}\mid c(\tilde x)\neq0, -h(\tilde x)\leq K\}$
is finite for all $K\in\mathbb R$. This is a finite dimensional
vector space over $\Lambda_\omega$, independent of the choice of $h$.
Note that for a section $s:\mathcal X\to\tilde{\mathcal X}$ the indicator
functions for $s(x)$, $x\in\mathcal X$, define a basis of $C^*(X;\Lambda_\omega)$.

To describe the differential let us call two elements 
$\gamma_1, \gamma_2\in\mathcal P_{x,y}$ equivalent if $\gamma_2^{-1}\gamma_1$
vanishes in $H_1(M;\mathbb R)$. Let $p_{x,y}:\mathcal P_{x,y}\to\mathcal
P_{x,y}'$ denote the projection onto the space of equivalence classes.
$\Gamma$ acts free and transitively on $\mathcal P_{x,y}'$.
The differential on $C^*(X;\Lambda_\omega)$ is determined by the counting functions
\begin{equation}\label{E:count}
\mathbb I_{x,y}':=(p_{x,y})_*\mathbb I_{x,y}:\mathcal P_{x,y}'\to\mathbb Z,
\qquad
\mathbb I_{x,y}'(a):=\sum_{p_{x,y}(\gamma)=a}\mathbb I_{x,y}(\gamma).
\end{equation}
Note that these sums are finite in view of Proposition~\ref{P:Nov}.

Now suppose $Y$ is a vector field on $W=M\times[-1,1]$ as in
section~\ref{SS:anom}. Assume $Y$ satisfies \MS\ and \LL.
Suppose $p^*\omega$ is Lyapunov for $Y$ where $p:W\to M$ denotes
the projection. As in section~\ref{SS:anom}
the differential of the Novikov complex $C^*(Y;\Lambda_{p^*\omega})$
gives rise to a homomorphism of Novikov complexes
\begin{equation}\label{E:uY}
u^Y:C^*(X_-;\Lambda_\omega)\to C^*(X_+;\Lambda_\omega).
\end{equation}
It is well known that \eqref{E:uY} is a quasi isomorphism.
Let $s_\pm:\mathcal X_\pm\to\tilde{\mathcal X}_\pm$ be sections
and equip $C^*(X_\pm;\Lambda_\omega)$ with the corresponding base.
Assume $(X_-,s_-)$ and $(X_+,s_+)$ determine the same Euler structure
\cite{T90, BH03}. Recall that this implies
\begin{equation}\label{E:eul}
\sum_{x\in\mathcal X_+}h^\eta(s(x))
-\sum_{x\in\mathcal X_-}h^\eta(s(x))
=\eta(\cs(X_-,X_+))
\end{equation}
for all $\eta\in\mathcal Z^1(M;\mathbb C)$ and all
smooth functions $h^\eta:\tilde M\to\mathbb C$ with $dh^\eta=\pi^*\eta$.

A result of Hutchings--Lee \cite{HL99} and Pajitnov \cite{P03}
tells that if $Y$ in addition satisfies \NCT, then the torsion of \eqref{E:uY} is
\begin{equation}\label{E:HLP}
\pm T(u^Y)=\pm\exp(h_*\mathbb P^{X_+}-h_*\mathbb P^{X_-})
\in 1+\Lambda^+_\omega.
\end{equation}

\subsection{Two lemmas}\label{SS:2lemmas}

Let $\Gamma:=\img\bigl(\pi_1(M)\to H_1(M;\mathbb R)\bigr)$, let
$\omega\in\mathcal Z^1(M;\mathbb R)$ be a closed one form,
suppose $\omega:\Gamma\to\mathbb R$ is injective and
let $\Lambda_\omega$ denote the Novikov field as introduced in
section~\ref{SS:HLP}.
For a closed one form $\eta\in\Omega^1(M;\mathbb C)$ we let $L^1_\eta$ 
denote the Banach algebra of all functions $\lambda:\Gamma\to\mathbb C$ with
$\|\lambda\|_\eta:=\sum_{\gamma\in\Gamma}|\lambda(\gamma)e^{\eta(\gamma)}|<\infty$
equipped with the convolution product. 
Moreover, let us write $\ev_\eta:L^1_\eta\to\mathbb C$ for the
homomorphism given by
$\ev_\eta(\lambda):=L(\lambda)(\eta)=\sum_{\gamma\in\Gamma}\lambda(\gamma)e^{\eta(\gamma)}$.

\begin{lemma}\label{L:invnov}
Suppose $0\neq\lambda\in\Lambda_\omega\cap L^1_\eta$. Then there exists
$t_0\in\mathbb R$ so that $\lambda^{-1}\in\Lambda_\omega\cap
L^1_{\eta+t\omega}$ for all $t\geq t_0$.
\end{lemma}

\begin{proof}
Using the Novikov property of $\lambda$ and the injectivity of
$\omega:\Gamma\to\mathbb R$ it is easy to see that we may w.l.o.g.\
assume $1-\lambda\in\Lambda_\omega^+$. Since $\lambda\in L^1_\eta$
we have $\|\lambda\|_\eta<\infty$. Using the Novikov property of $\lambda$
and the injectivity of $\omega:\Gamma\to\mathbb R$ again we find
$t_0\in\mathbb R$ so that $\|1-\lambda\|_{\eta+t\omega}<1$, for all $t\geq t_0$.
Since $L^1_{\eta+t\omega}$ is a Banach algebra 
$\sum_{k\geq0}(1-\lambda)^k$ will converge and 
$\lambda^{-1}\in\Lambda_\omega\cap L^1_{\eta+t\omega}$, for all $t\geq t_0$.
\end{proof}

Recall that we have a bijection $\exp:\Lambda_\omega^+\to1+\Lambda_\omega^+$.

\begin{lemma}\label{L:log}
Suppose $\lambda\in\Lambda_\omega^+$ and $\exp(\lambda)\in L_\eta^1$.
Then there exists $t_0\in\mathbb R$ so that $\lambda\in
\Lambda_\omega\cap L^1_{\eta+t\omega}$, for all $t\geq t_0$.
\end{lemma}

\begin{proof}
Similar to the proof of Lemma~\ref{L:invnov} using 
$\log(1-\mu)=-\sum_{k>0}\frac{\mu^k}k$.
\end{proof}

\subsection{Computation of the anomaly}\label{SS:can}

With the help of the Hutchings--Lee formula it is possible
to compute the right hand side of \eqref{E:anom} in terms of closed
trajectories under some assumptions.

\begin{proposition}\label{P:inp_anom}
Suppose $Y$ is a vector field on $M\times[-1,1]$ as in
Proposition~\ref{P:anom} which satisfies \MS, \LL, \NCT\ and \EG.
Let $\eta\in\mathcal Z^1(M;\mathbb C)$ be a closed one form.
Suppose $\omega\in\mathcal Z^1(M;\mathbb R)$ such that
$\omega:\Gamma\to\mathbb R$ is injective
and such that $[p^*\omega]$ is a Lyapunov class for $Y$.
Then there exists $t_0$ such that for $t>t_0$ 
we have
$\eta+t\omega\in\bigl(\mathring\RX^{X_+}\setminus\Sigma^{X_+}\bigr)
\cap\bigl(\mathring\RX^{X_-}\setminus\Sigma^{X_-}\bigr)$,
$L(h_*\mathbb P^{X_+}-h_*\mathbb P^{X_-})(\eta+t\omega)$
converges absolutely, and
$$
\frac{(T\Int_{\eta+t\omega}^{X_+})^2}{(T\Int_{\eta+t\omega}^{X_-})^2}
=\bigl(e^{L(h_*\mathbb P^{X_+}-h_*\mathbb P^{X_-})(\eta+t\omega)}\bigr)^2.
$$
\end{proposition}

\begin{proof}
Since $X_\pm$ satisfies \EG, and since the cohomology class of $\omega$
contains a Lyapunov form for $X_\pm$ we obtain from
Proposition~\ref{P:expR}, Theorem~\ref{T:int} and Proposition~\ref{P:basic:R} that
$\eta+t\omega\in\bigl(\mathring\RX^{X_+}\setminus\Sigma^{X_+}\bigr)
\cap\bigl(\mathring\RX^{X_-}\setminus\Sigma^{X_-}\bigr)$ for sufficiently large $t$.
Arguing similarly for $Y$ we see that $p^*(\eta+t\omega)\in\RX^Y$
for sufficiently large $t$. 
Particularly, Proposition~\ref{P:anom} is applicable and we get
$$
\frac{(T\Int_{\eta+t\omega}^{X_+})^2}{(T\Int_{\eta+t\omega}^{X_-})^2}
=(Tu^Y_{\eta+t\omega})^2\cdot\bigl(e^{-(\eta+t\omega)(\cs(X_-,X_+))}\bigr)^2.
$$

Since $\RX^Y\subseteq\QX^Y$, see Theorem~\ref{T:int}, 
the Novikov complex of $Y$ is defined over the ring $\Lambda_t:=
\Lambda_\omega\cap L^1_{\eta+t\omega}$ for sufficiently large $t$. 
More precisely, for sufficiently large $t$
the counting functions \eqref{E:count} actually define a 
complex $C^*(Y;\Lambda_t)$ over $\Lambda_t$ with
$$
C^*(Y;\Lambda_\omega)
=C^*(Y;\Lambda_t)\otimes_{\Lambda_t}\Lambda_\omega.
$$
Since the basis determined by sections $s_\pm:\mathcal X_\pm\to\tilde{\mathcal X}_\pm$
obviously consist of elements in $C^*(Y;\Lambda_t)$ we
conclude that the torsion $T(u^Y)$ is contained in the quotient field
$Q(\Lambda_t)\subseteq\Lambda_\omega$. In view of \eqref{E:HLP}
Lemma~\ref{L:invnov} and Lemma~\ref{L:log} we thus have
$$
h_*\mathbb P^{X_+}-h_*\mathbb P^{X_-}\in\Lambda_t,
$$
and hence $L(h_*\mathbb P^{X_+}-h_*\mathbb P^{X_-})(\eta+t\omega)$
converges absolutely for sufficiently large $t$.

For sufficiently large $t$, 
let us write $\ev_t:\Lambda_t\to\mathbb C$ for the homomorphism given by
$\ev_t(\lambda):=L(\lambda)(\eta+t\omega)
=\sum_{\gamma\in\Gamma}\lambda(\gamma)e^{(\eta+t\omega)(\gamma)}$.
Clearly,
$$
C^*_{\eta+t\omega}(Y;\mathbb C)
=C^*(Y;\Lambda_t)\otimes_{\ev_t}\mathbb C.
$$
Moreover, using \eqref{E:eul} and Lemma~\ref{L:invnov}, it is easy to see that
this implies 
$$
\pm Tu^Y_{\eta+t\omega}\cdot e^{-(\eta+t\omega)(\cs(X_-,X_+))}
=\pm L(Tu^Y)(\eta+t\omega),
$$
and \eqref{E:HLP} yields
\[
\pm Tu^Y_{\eta+t\omega}\cdot e^{-(\eta+t\omega)(\cs(X_-,X_+))}
=\pm e^{L(h_*\mathbb P^{X_+}-h_*\mathbb P^{X_-})(\eta+t\omega)}.
\qedhere
\]
\end{proof}

\subsection{Proof of Theorem~\ref{T:tor}}\label{SS:Pthmint}

Let $X$ be a vector field satisfying \SEG.
Choose a vector field $Y$ on $M\times[-1,1]$ as in Definition~\ref{D:SEG}.
Note that $h_*\mathbb P^{X_+}=h_*\mathbb P^X$, $h_*\mathbb P^{X_-}=0$,
and $(T\Int_{\tilde\eta}^{X_-})^2=1$ for all $\tilde\eta\in\mathcal
Z^1(M;\mathbb C)$, see \eqref{E:BZ}.
Let $\Gamma:=\img\bigl(\pi_1(M)\to H_1(M;\mathbb R)\bigr)$.
Let $\omega_0$ be a Lyapunov form for $X$ such that
$\omega_0:\Gamma\to\mathbb R$ is injective, see 
Proposition~\ref{P:cone}.
>From Propositions~\ref{P:hl} (Appendix B) and \ref{P:cone}
we see that $[p^*\omega_0]$ is a Lyapunov class for $Y$.

Let $\eta\in\mathcal Z^1(M;\mathbb C)$ be a closed one form.
In view of Proposition~\ref{P:inp_anom} we have $\eta+t\omega_0\in\mathring\RX^X$
for sufficiently large $t$. Let $\omega$ be an arbitrary Lyapunov form for
$X$. Using Propositions~\ref{P:cone} and \ref{P:basic:R} we conclude that
$\eta+t\omega\in\mathring\RX^X$ for sufficiently large $t$.
Applying Proposition~\ref{P:inp_anom} to various $\eta$ we see that
\begin{equation}\label{E:777}
(T\Int^X_{\tilde\eta})^2
=(e^{L(h_*\mathbb P^X)(\tilde\eta)})^2
\end{equation}
holds for an open subset of
$\tilde\eta\in(\mathring\RX^X\setminus\Sigma^X)\cap\mathring\PX^X$.
By analyticity, see Propositions~\ref{P:basic:P} and \ref{P:Tint}, equality
\eqref{E:777} holds for all
$\tilde\eta\in(\mathring\RX^X\setminus\Sigma^X)\cap\mathring\PX^X$.
Using \eqref{E:cont:P} and \eqref{E:cont:int}
we see that the relation \eqref{E:777}
continues to hold for all $\tilde\eta\in(\RX^X\setminus\Sigma^X)\cap\PX^X$.
This completes the proof of Theorem~\ref{T:tor}.

\subsection{Proof of Theorem~\ref{T:tor}'}\label{SS:Pthmintp}

Let $\omega$ be a Lyapunov form for $X$ and assume $\omega:\Gamma\to\mathbb
R$ is injective, see Proposition~\ref{P:cone}. Let $\eta\in\mathcal
Z^1(M;\mathbb C)$. Note that since $H_{\eta_0}(M;\mathbb C)=0$ the
deRham cohomology will be acyclic, generically. More precisely, 
$H^*_{\eta+t\omega}(M;\mathbb C)=0$ for sufficiently large $t$.
In view of Theorem~\ref{T:int} we also have $H^*_{\eta+t\omega}(X;\mathbb
C)=0$ for sufficiently large $t$.

As in section~\ref{SS:can} one shows that for sufficiently large $t$
the Novikov complex $C^*(M;\Lambda_\omega)$
is actually defined over the ring $\Lambda_t:=\Lambda_\omega\cap
L^1_{\eta+t\omega}$,
$$
C^*(X;\Lambda_\omega)=C^*(X;\Lambda_t)\otimes_{\Lambda_t}\Lambda_\omega.
$$
Moreover, for sufficiently large $t$
$$
C^*(X;\Lambda_t)\otimes_{\ev_t}\mathbb C=C^*_{\eta+t\omega}(X;\mathbb C).
$$
We conclude that the Novikov complex $C^*(X;\Lambda_\omega)$ is acyclic.

Let $Y=-\grad_{g_0}f$ be a Morse--Smale vector field with zero set
$\mathcal Y$. Let $s_X:\mathcal X\to\tilde{\mathcal X}$ and
$s_Y:\mathcal Y\to\tilde{\mathcal Y}$ be a sections and assume that they
define the same Euler structure, i.e.\
$$
\sum_{x\in\mathcal X}h^{\tilde\eta}(s_X(x))
-\sum_{y\in\mathcal Y}h^{\tilde\eta}(s_Y(y))
=\tilde\eta(\cs(X,Y))
$$
for all $\tilde\eta\in\mathcal Z^1(M;\mathbb C)$ and all smooth functions
$h:\tilde M\to\mathbb C$ with $dh^{\tilde\eta}=\pi^*\tilde\eta$, see section~\ref{SS:HLP}.
Equip the complexes
$C^*(X;\Lambda_\omega)$ and $C^*(Y;\Lambda_\omega)$ with the corresponding graded
bases. For the torsion we have \cite{HL99, P03}
$$
\frac{\pm T\bigl(C^*(X;\Lambda_\omega)\bigr)}
{\pm T\bigl(C^*(Y;\Lambda_\omega)\bigr)}
=\pm\exp(h_*\mathbb P^X)\in1+\Lambda^+_\omega.
$$
As in section~\ref{SS:can} one shows that this torsion must be contained in the quotient field
$Q(\Lambda_t)\subseteq\Lambda_\omega$, hence $(h_*\mathbb P)(\eta+t\omega)$
converges absolutely, and thus
$\eta+t\omega\in\mathring\RX\setminus\Sigma$ for sufficiently large $t$.
Again, this remains true for arbitrary Lyapunov $\omega$ in view of
Proposition~\ref{P:cone}.

Equip the complexes $C^*_{\eta+t\omega}(X;\mathbb C)$ and
$C^*_{\eta+t\omega}(Y;\mathbb C)$ with the graded bases determined by the
indicator functions. As in section~\ref{SS:can} we conclude that
$$
\frac{\pm T\bigl(C^*_{\eta+t\omega}(X;\mathbb C)\bigr)}
{\pm T\bigl(C^*_{\eta+t\omega}(Y;\mathbb C)\bigr)}
e^{-(\eta+t\omega)(\cs(X,Y))}
=\pm e^{L(h_*\mathbb P^X)(\eta+t\omega)}
$$
for sufficiently large $t$. Using \eqref{EE:38} this implies
$$
\frac{(T\Int_{\eta+t\omega}^X)^2}
{(T\Int_{\eta+t\omega}^Y)^2}
=(e^{L(h_*\mathbb P^X)(\eta+t\omega)})^2.
$$
In view of \eqref{E:BZ} we have $(T\Int_{\eta+t\omega}^Y)^2=1$. 
We conclude that 
$$
(T\Int^X_{\tilde\eta})^2=(e^{L(h_*\mathbb P^X)(\tilde\eta)})^2
$$
holds for an open set of
$\tilde\eta\in(\mathring\RX^X\setminus\Sigma^X)\cap\mathring\PX^X$.
By analyticity, see Propositions~\ref{P:basic:P} and \ref{P:Tint}, this
relation holds for
all $\tilde\eta\in(\mathring\RX^X\setminus\Sigma^X)\cap\mathring\PX^X$.
Using \eqref{E:cont:int} and \eqref{E:cont:P} it remains true for all
$\eta\in(\RX^X\setminus\Sigma^X)\cap\PX^X$.

\begin{appendix}

\section{Proof of Proposition~\ref{P:cano}}\label{S:cano}

We will make use of the following lemma whose proof we leave to the reader.

\begin{lemma}\label{L:extpF}
Let $N$ be a compact smooth manifold, possibly with boundary, and let
$K\subseteq N$ be a compact subset. Let 
$L:=N\times\partial I\cup K\times I$ where $I:=[0,1]$. Suppose
$F$ is a smooth function defined in a neighborhood of $L$ so that
$\partial F/\partial t<0$ whenever defined, and so that
$F(x,0)>F(x,1)$ for all $x\in N$. Then there exists a smooth function
$G:N\times I\to\mathbb R$ which agrees with $F$ on a neighborhood of $L$
and satisfies $\partial G/\partial t<0$.
\end{lemma}

For $\rho\geq\epsilon>0$ define
$$
\mathbb D_{\rho,\epsilon}:=\bigl\{
(y,z)\in\mathbb R^q\times\mathbb R^{n-q}\bigm|
-\rho\leq-\tfrac12|y|^2+\tfrac12|z|^2\leq\rho, |y|\cdot|z|\leq\epsilon
\bigr\}.
$$

\begin{lemma}\label{L:extF}
Suppose $F:\mathbb D_{\rho,\rho}\to\mathbb R$ is a smooth function 
with $F(0)=0$ which is strictly decreasing along non-constant 
trajectories of $X$, see \eqref{E:4}.
Then there exists $\rho>\epsilon>0$ and a smooth function
$G:\mathbb D_{\rho,\rho}\to\mathbb R$ which is strictly decreasing along non-constant
trajectories of $X$, which coincides with $F$ on a neighborhood of 
$\partial\mathbb D_{\rho,\rho}$ and which coincides with 
$-\frac12|y|^2+\frac12|z|^2$ on $\mathbb D_{\epsilon,\epsilon}$.
\end{lemma}

\begin{proof}
Consider the partially defined function which coincides with 
$F$ in a neighborhood 
of $\partial\mathbb D_{\rho,\rho}$ and which coincides with
$-\frac12|y|^2+\frac12|z|^2$ on a neighborhood of
$\mathbb D_{\epsilon,\epsilon}$.
We will extend this to a globally defined smooth function
$\mathbb D_{\rho,\rho}\to\mathbb R$ which is strictly decreasing along non-constant
trajectories of $X$. This will be accomplished in two steps.

For the first step notice
that $\overline{\mathbb D_{\rho,\epsilon}\setminus\mathbb D_{\epsilon,\epsilon}}$
is diffeomorphic to $N\times I$ where 
$N=S^{q-1}\times D^{n-q}\cup D^q\times S^{n-q-1}$. Here $S^{k-1}$ and $D^k$ denote
unit sphere and unite ball in $\mathbb R^k$, respectively.
Choosing $\epsilon$ sufficiently small we can apply Lemma~\ref{L:extpF}
with $K=\emptyset$, and obtain an extension to $\mathbb D_{\rho,\epsilon}$.

For the second step notice that 
$\overline{\mathbb D_{\rho,\rho}\setminus\mathbb D_{\rho,\epsilon}}$
is diffeomorphic to $N\times I$ where $N=C^q\times S^{n-q-1}$
and $C^q:=\{y\in\mathbb R^q\mid 1\leq|y|\leq2\}$. Applying
Lemma~\ref{L:extpF} with $K=\partial C^q$, provides the desired extension
to $\mathbb D_{\rho,\rho}$.
\end{proof}

\begin{proof}[Proof of Proposition~\ref{P:cano}]
Let $\omega\in\Omega^1(M;\mathbb R)$ be a closed one form such that
$\omega(X)<0$ on $M\setminus\mathcal X$. 
By adding a small closed one form with support
contained in $M\setminus\mathcal X$ we may in addition assume that
the cohomology class of $\omega$ is rational. Multiplying with
a positive number we may assume that the cohomology class of $\omega$ is integral.
Moreover, in view of Lemma~\ref{L:extF}
we may assume that $\omega$
has canonical form in a neighborhood of $\mathcal X$. More precisely,
for every $x\in\mathcal X_q$ there exist coordinates $(x_1,\dotsc,x_n)$
centered at $x$ in which 
\begin{equation}\label{E:98}
X=\sum_{i\leq q}x_i\frac\partial{\partial x_i}
-\sum_{i>q}x_i\frac\partial{\partial x_i}
\qquad\text{and}\qquad
\omega=-\sum_{i\leq q}x_idx^i+\sum_{i>q}x_idx^i.
\end{equation}
Define a Riemannian metric $g$ on $M$ as follows. 
On a neighborhood of $\mathcal X$ on which $X$ and
$\omega$ have canonic form define  $g:=\sum_i(dx^i)^2$. Note that this
implies $\omega=-g(X,\cdot)$ where defined. Since $\omega(X)<0$ we have 
$TM=\ker\omega\oplus[X]$ over $M\setminus\mathcal X$. Extend 
$g|_{\ker\omega}$ smoothly to a fiber metric on $\ker\omega$ over
$M\setminus\mathcal X$, and let the restriction of $g$ to $\ker\omega$ 
be given by this extension. Moreover, set $g(X,X):=-\omega(X)$ and
$g(X,\ker\omega):=0$. This defines a smooth Riemannian metric on $M$,
and certainly $\omega=-g(X,\cdot)$.
\end{proof}

\section{Vector fields on $M\times[-1,1]$}\label{app:B}

\begin{proposition}\label{P:hl}
Let $X_\pm$ be two vector fields on $M$.
Then there exists a vector field $Y$ on $M\times[-1,1]$
such that $Y(z,s)=X_+(z)+(s-1)\partial/\partial s$
in a neighborhood of $\partial_+W$, such that
$Y(z,s)=X_-(z)+(-s-1)\partial/\partial s$ in a neighborhood of
$\partial_-W$, and such that $ds(Y)<0$ on $M\times(-1,1)$.
Moreover, every such vector field has the following property:
If $\xi\in H^1(M;\mathbb R)$ is a Lyapunov class for $X_+$ and $X_-$,
then $p^*\xi\in H^1(M\times[-1,1];\mathbb R)$ is a Lyapunov class for $Y$,
where $p:M\times[-1,1]\to M$ denotes the projection.
\end{proposition}

\begin{proof}
The existence of such a vector field $Y$ is obvious.
Suppose $\xi\in H^1(M;\mathbb R)$ is a Lyapunov class for $X_-$ and $X_+$.
It is easy to construct a closed one form $\omega\in\mathcal
Z^1(M\times[-1,1];\mathbb R)$ representing $p^*\xi$ such that
$\omega_\pm(X_\pm)<0$ on $M\setminus\mathcal X_\pm$, where 
$\omega_\pm:=\iota_\pm^*\omega\in\mathcal Z^1(M;\mathbb R)$, and
$\iota_\pm:M\to M\times\{\pm1\}\subseteq M\times[-1,1]$ denotes
the canonic inclusions.
We may moreover assume that $i_{\partial_s}\omega$ vanishes in a
neighborhood of $M\times\{\pm1\}$.
For sufficiently large $t$ the form $\omega+tds\in\mathcal
Z^1(M\times[-1,1];\mathbb R)$ will be a Lyapunov form for $Y$
representing $p^*\xi$.
\end{proof}

\end{appendix}

\end{document}